\documentclass[preprint]{elsarticle}
\usepackage{lineno}

\journal{JOUARNAL's NAME}

\usepackage{slashbox}
\usepackage{amsmath}
\usepackage{amssymb}
\usepackage{amsthm}
\usepackage{graphicx}
\usepackage{epic,eepic,epsfig}
\usepackage{color}
\usepackage{subfigure}  
\usepackage{placeins}   
\usepackage{multirow}
\usepackage{epstopdf}
\usepackage{varwidth}
\usepackage{float}
\usepackage[table]{xcolor}
\usepackage{bm}
\usepackage[title]{appendix}
\usepackage{diagbox}
\usepackage{booktabs}

\usepackage{algorithm}
\usepackage{algpseudocode}

\newtheorem{theorem}{Theorem}[section]
\newtheorem{example}{Example}[section]
\newtheorem{lemma}{Lemma}[section]
\newtheorem{remark}{Remark}[section]

\newtheorem{assumption}{Assumption}[section]

\definecolor{tabclr}{cmyk}{0,0,1,0}

\allowdisplaybreaks[4]

\usepackage{algorithm}
\usepackage{algpseudocode}

\allowdisplaybreaks[4]

\makeatletter
\def\ps@pprintTitle{%
 \let\@oddhead\@empty
 \let\@evenhead\@empty
 \let\@oddfoot\@empty
 \let\@evenfoot\@empty
}
\makeatother

\begin{document}
\title{Exponential Runge--Kutta methods for parabolic equations with state-dependent delay}

\author[bjut]{Qiumei Huang}
\ead{qmhuang@bjut.edu.cn}

\author[ins]{Alexander Ostermann}
\ead{alexander.ostermann@uibk.ac.at}

\author[bjut,ins]{Gangfan Zhong}
\ead{gfzhong@emails.bjut.edu.cn}

\address[bjut]{School of Mathematics, Statistics and Mechanics, Beijing University of Technology, 100124 Beijing, China}
 
\address[ins]{Department of Mathematics, University of Innsbruck, Technikerstr.~13, 6020 Innsbruck, Austria}


\begin{abstract}
The aim of this paper is to construct and analyze exponential Runge--Kutta methods for the temporal discretization of a class of semilinear parabolic problems with arbitrary state-dependent delay. First, the well-posedness of the problem is established. Subsequently, first and second order schemes are constructed. They are based on the explicit exponential Runge--Kutta methods, where the delayed solution is approximated by a continuous extension of the time discrete solution. Schemes of arbitrary order can be constructed using the methods of collocation type. The unique solvability and convergence of the proposed schemes are established. Finally, we discuss implementation issues and present some numerical experiments to illustrate our theoretical results.
\end{abstract}


\maketitle
%
%

\section{Introduction}\label{Sec:intro}
\setcounter{equation}{0}

Delay differential equations (DDEs) and delay partial differential equations are important tools in modeling real-world processes with inherent time delays, including problems in physics, chemistry, control theory, biology, and other fields. Compared with non-delay models, delay models generally provide a more realistic description of the dynamic nature of real-world systems. 
For simplicity, delays are often assumed to be constant. However, this assumption rarely applies to systems in practice, where delays can be time- or even state-dependent. For a detailed overview of state-dependent DDEs, we refer the reader to \cite{HartungKrisztin06Book}. The numerical analysis for state-dependent DDEs has also been well-developed; see \cite{BellenMaset09:1,BellenZennaro13Book}. 

In contrast, the study of partial differential equations with state-dependent delay remains an active area of research, focusing on the theory of dynamical systems \cite{HernandezFernandes21:753,HernandezPierri16:6856,KrisztinRezounenko16:4454,Rezounenko09:3978}.
In this study, we consider the numerical solution of the following class of (abstract) semilinear parabolic problems with state-dependent delay
\begin{equation}\label{Eqn:problem}
\left\{\begin{aligned}
&u^\prime(t) + Au(t) = f\big(t,u(t),u \big(t-\tau(t,u(t))\big) \big), && t > 0, \\
&u(t) = \phi(t)  ,   && t\leq 0,
\end{aligned}\right.
\end{equation}
where the delay $\tau$ depends on time and the actual state. In the past, significant numerical research has been conducted for problems of the form \eqref{Eqn:problem} in the case of constant delay; see, for example, \cite{HuangVandewalle12:579,vanPieter86:1,XuHuang23:57}.  Operator splitting for abstract delay equations has also been investigated in \cite{BatkaiCsomos13:315,CsomosNickel08:2234,HansenStillfjord14:673}. 

In recent years, exponential integrators have attracted considerable attention due to their effectiveness for stiff semilinear systems.  By treating the linear term exactly and approximating the nonlinearity in an explicitly way, they are able to solve stiff problems in an accurate and efficient way.  For a comprehensive overview of exponential integrators, we refer the reader to \cite{HochbruckOstermann10:209}. Stability and convergence of exponential integrators for DDEs with constant delay have been studied in \cite{XuZhao11:1089,ZhaoZhan16:96,ZhaoZhan18:45}. 
Using the sun-star theory, And\`o and Vermiglio \cite{AndoVermiglio25:1842} reformulated DDEs as abstract ordinary differential equations, making them amenable to exponential Runge--Kutta (ERK) methods.
The ERK methods have
also been applied to semilinear parabolic problems with constant delay \cite{DaiHuang23:350,Zhan21:113279} and (non-vanishing) time-dependent delay \cite{HuangOstermann25}. 

In this paper, we aim to extend the ERK methods as presented in \cite{HochbruckOstermann05:1069,HochbruckOstermann05:323} to the larger class of problems \eqref{Eqn:problem}. The core idea is to construct continuous extensions of the discrete solutions obtained by ERK methods to approximate the delayed solution. 
As far as we are aware, this is the first study to address the numerical analysis of partial differential equations with state-dependent delay.

The outline of the paper is as follows. 
In Section~\ref{Sec:framework}, we summarize the employed abstract framework and establish the well-posedness of the initial value problem \eqref{Eqn:problem}. Further, we construct the exponential Euler method for \eqref{Eqn:problem} and analyze the convergence in Section~\ref{Sec:Euler}. Based on an explicit ERK method, a second order method for \eqref{Eqn:problem} is presented in Section~\ref{Sec:second}. In Section~\ref{Sec:collocation}, we construct $s$-stage ERK methods of collocation type for \eqref{Eqn:problem} and establish their unique solvability. It is shown that the methods achieve order $s$ and can further achieve superconvergence provided that underlying quadrature rule is of order $s+1$.   Finally, we discuss the implementation of the proposed methods and present some numerical experiments in Section~\ref{Sec:experiments} to illustrate the theoretical results.

\section{Analytical Framework}\label{Sec:framework}

Our analysis below will be based on an abstract formulation of \eqref{Eqn:problem} as an evolution equation with delay in a Banach space $(X,\|\cdot\|)$. Let $D(A)$ denote the domain of $A$ in $X$. Our basic assumptions on the operator are as follows.

\begin{assumption}\label{Ass:sectorial}
Let the operator $A: D(A) \rightarrow X$ be an infinitesimal generator of  a compact analytic semigroup $\mathrm{e}^{-tA}$ in $X$,  and let $D(A)$ be dense in $X$. Without restriction of generality, we assume that the spectrum $\sigma(A)$ of $A$, satisfies $\mathrm{Re}\,\sigma(A)>0$.
\end{assumption}

Under this assumption, the fractional powers of $A$ are well defined. We recall that $A$ satisfies the properties (see \cite{Henry81Book,Lunardi95Book})
$$
\|\mathrm{e}^{-t A}\|_{X \leftarrow X}+ \|t^\alpha A^\alpha \mathrm{e}^{-t A} \|_{X \leftarrow X}   \leq C, \quad  \alpha,t \geq 0.
$$
It follows that the $\varphi_k$ functions appearing in exponential integrators, defined by
$$
\varphi_k(-tA)=\frac{1}{t^k}\int_0^t \mathrm{e}^{-(t-\xi)A}\frac{\xi^{k-1}}{(k-1)!}\,\mathrm{d}\xi,\quad k\geq1,
$$
satisfy $\|\varphi_k(-tA)\|\leq C$ for $t\geq 0$.

Our basic assumptions on $f$ and $\tau$ are stated below. 
\begin{assumption}\label{Ass:lip}
Let the nonlinearity $f:[0, +\infty) \times X \times X \rightarrow X$ and the delay function $\tau:[0,+\infty)\times X\to[0,+\infty)$ be Lipschitz continuous.
\end{assumption}

This assumption infers that there exist real numbers $L_f$ and $L_\tau$ such that
$$
\begin{aligned}
\|f(t_1, v_1,w_1)-f(t_2,v_2, w_2)\| & \leq L_f  ( |t_1-t_2|+\|v_1-v_2\|+\|w_1-w_2\| ), \\
|\tau(t_1,v_1)-\tau(t_2,v_2)|& \leq L_\tau (|t_1-t_2|+ \|v_1-v_2\|),
\end{aligned}
$$
for all $t_1,t_2 \in[0, +\infty)$ and all $v_1,v_2,w_1,w_2\in X$.

Given an interval $\mathcal{I}$, we denote $C_{U}(\mathcal{I};X)$ as the space of uniformly continuous functions on $\mathcal{I}$ equipped the supremum norm. The H\"older spaces $C^\alpha(\mathcal{I};X)$ $(0<\alpha\le1)$ 
and the Lipschitz spaces $C^{k,1}(\mathcal{I};X)$ ($k\in \{0\}\cup\mathbb{N}$) are defined in the usual way, and their 
norms are denoted by 
$\|\cdot\|_{C^\alpha(\mathcal{I};X)}$ and 
$\|\cdot\|_{C^{k,1}(\mathcal{I};X)}$, respectively. 
For convenience, we recall that 
$$
\|u\|_{C^\alpha(\mathcal{I};X)}
=\sup_{t\in\mathcal{I}}\|u(t)\|+[u]_{C^\alpha(\mathcal{I};X)},\quad \|u\|_{C^{k,1}(\mathcal{I};X)}=\sum_{|\beta|\le k}\|\partial^\beta u\|_{C(\mathcal{I};X)}+[\partial^k u]_{C_{Lip}(\mathcal{I};X)},
$$ 
where 
$$
[u]_{C^\alpha(\mathcal{I};X)}
=\sup_{s,t\in\mathcal{I},s\ne t}\frac{\|u(t)-u(s)\|}{|t-s|^\alpha},\quad
[u]_{C_{Lip}(\mathcal{I};X)}=\sup_{s,t\in\mathcal{I},s\ne t}\frac{\|u(t)-u(s)\|}{|t-s|}.
$$
The existence and uniqueness of solutions to the initial value problem \eqref{Eqn:problem} is given by the following theorem. 

\begin{theorem}\label{Thm:well-posedness}
Under Assumptions \ref{Ass:sectorial}-\ref{Ass:lip}, if $\phi(t)$ is Lipschitz continuous for $t\leq 0$ and $\phi(0)\in D(A)$, then there exists a time $T=T(\phi)>0$ such that  initial value problem \eqref{Eqn:problem} has a unique solution $u \in C_{U}( (-\infty,T];X)\cap C^1([0,T];X)$.
\end{theorem}
\begin{proof}
We denote by $u_t$ the element of $C_{U}((-\infty,0];X)$ defined by the formula $u_t(\theta)= u(t+\theta)$ for $\theta\in (-\infty,0]$.
Let $F:[0,+\infty)\times C_{U}((-\infty,0];X) \to X$ be the function defined by $F(t,u_t)=f\big(t,u_t(0),u_t(-\tau(t,u_t(0)) ) \big)$. Then the problem \eqref{Eqn:problem} can be reformulated as
$$
u^\prime(t)+Au(t)=F(t,u_t),\quad u_0=u|_{(-\infty,0]}=\phi\in C_{U}((-\infty,0];X).
$$
Since $F$ is continuous, it follows from \cite{Fitzgibbon78:1} that there exists a positive time $T=T(\phi)$ and a function $u \in C_{U}( (-\infty,T];X)$ such that
\begin{equation}\label{Eqn:mild}
u(t)=\mathrm{e}^{-tA}\phi(0) + \int_0^t \mathrm{e}^{-(t-s)A}F(s,u_s)\,\mathrm{d}s,\quad t\in [0,T],
\end{equation}
which is a so called {mild solution}. As the initial function $\phi(t)$ is Lipschitz continuous for $t\leq 0$ and $\phi(0)\in D(A)$, by following the idea in \cite[Section 2]{KrisztinRezounenko16:4454} one can establish the uniqueness of the solution as below. 

Noting that $g(t)=F(t,u_t)$ belongs to $C([0,T];X)$, the initial value problem
$$
v^\prime(t)+Av(t)=g(t),\quad v(0)=\phi(0)
$$ 
admits a unique mild solution $v=u$. Since $u(0)\in D(A) \subset  D(A^\frac{1}{2})$,
it follows from \cite[Corollary 4.2.2]{Lunardi95Book} that $u\in C^{\frac{1}{2}}([0,T];X)$. Noting the fact that $\phi$ is Lipschitz continuous for $t\leq 0$ and using the Lipschitz conditions of $f$ and $\tau$, one has, for $0\leq s < t \leq T$, 
$$
\begin{aligned}
\|g(t)-g(s)\| &\leq L_f\big(|t-s|+\|u(t)-u(s)\|+\big\|u\big(t-\tau(t,u(t)) \big) - u\big(s-\tau(t,v(s)) \big) \big\| \big)\\
&\leq L_f\big(|t-s|+ [u]_{C^{\frac{1}{2}}}  |t-s|^{\frac{1}{2}} + [u]_{C^{\frac{1}{2}}}|t-s -\tau(t,u(t))+\tau(s,u(s))|^{\frac{1}{2}} \big),
\end{aligned}
$$
where $[u]_{C^{\frac{1}{2}}}$ denotes the H\"older seminorm taken over $(-\infty,T]$.
By the relation $|t-s|\leq T^{\frac{1}{2}}|t-s|^{\frac{1}{2}}$, we have
$$
\begin{aligned}
|t-s -\tau(t,u(t))+\tau(s,u(s))|^{\frac{1}{2}} & \leq \big(  |t-s|  +L_{\tau}|t-s|+ L_{\tau}\|u(t)-u(s)\| \big)^{\frac{1}{2}} \\
& \leq \big( T^{\frac{1}{2}} + L_{\tau}T^{\frac{1}{2}} +L_{\tau}[u]_{C^{\frac{1}{2}}}\big)^{\frac{1}{2}}|t-s|^{\frac{1}{4}}.
\end{aligned}
$$
A combination of the above two inequalities yields $g\in C^{\frac{1}{4}}([0,T];X)$. By \cite[Theorem 4.3.1]{Lunardi95Book}, the mild solution $u$ is strict, i.e., $u\in  C^1([0,T];X)\cap C( [0,T];D(A))$. The uniqueness of the solution is addressed next. If $u$ and $w$ are two mild solutions, we have
$$
\begin{aligned}
&\big\|f \big(t,u(t),u\big(t-\tau(t,u(t)) \big) \big)-f \big(t,w(t),w\big(t-\tau(t,w(t)) \big) \big)\big\| \\
&\qquad \leq L_f\|u(t)-w(t)\|+L_f \big\|u\big(t-\tau(t,u(t))\big) - w\big(t-\tau(t,w(t)) \big)\big\| \\
&\qquad \leq L_f\|u(t)-w(t)\|+L_fL_uL_{\tau}\|u(t)-w(t)\| + L_f\big\|u\big(t-\tau(t,w(t))\big) - w\big(t-\tau(t,w(t)) \big)\big\|, 
\end{aligned}
$$
where $L_u$ is the Lipschitz constant of $u$ in $(-\infty,T]$. It follows that, for $t\in[0,T]$,
$$
\|f \big(t,u(t),u\big(t-\tau(t,u(t)) \big)-f \big(t,w(t),w\big(t-\tau(t,w(t)) \big)\big\| 
\leq  \big( 2L_f +L_fL_uL_{\tau} \big) \|u-w\|_{C_{U}((-\infty,t];X)}.
$$
Recalling \eqref{Eqn:mild}, the difference of the two solutions $u$ and $w$ is bounded by
$$
\|u-w\|_{C_{U}((-\infty,t];X)} \leq \big( 2L_f +L_fL_uL_{\tau} \big)\int_0^t \|\mathrm{e}^{-(t-s)A}\|_{X \leftarrow X}\|u-w\|_{C_{U}((-\infty,s];X)}\,\mathrm{d}s,\quad t\in[0,T].
$$
The uniqueness follows from the boundedness of the semigroup and Gronwall's  inequality.
\end{proof}

\section{Exponential Euler method}\label{Sec:Euler}
In this section, we employ the exponential Euler method for the initial value problem \eqref{Eqn:problem} and analyze its convergence.

Let $I_h=\{t_n:0=t_0<t_1<\cdots<t_N=T\}$ be a mesh for the time domain $[0,T]$, and set
$
h_{n+1}=t_{n+1}-t_n$ and $h=\max_{1\leq n \leq N} h_n.
$
We first construct the exponential Euler approximation of $u(t_1)$. For this purpose, we consider the following problem
\begin{equation}\label{Eqn:euler-w1}
\left\{\begin{aligned}
&w_1^\prime(t)+Aw_1(t)=g(t,w_1(t)),\quad  t\in[0,t_1],\\
&w_1(0)=\phi(0),\end{aligned}\right.
\end{equation}
where $g(t,w_1(t))=f\big(t,w_1(t),\psi \big(t-\tau(t,w_1(t))  \big) \big)$ with $\psi$ defined by
$$
\psi(t)=\left\{\begin{aligned} & \phi(t), &&t\in (-\infty,0], \\
&w_1(t), &&t\in[0,t_1]. \end{aligned}\right.
$$
Applying the exponential Euler method \cite{HochbruckOstermann05:1069} to \eqref{Eqn:euler-w1} gives the following approximation $u_1$ to $u(t_1)$: 
$$
u_1 = \mathrm{e}^{-  h_1A}\phi(0) +   h_1\varphi_1(-h_1A)f\big(0,\phi(0),\phi (-\tau(0,\phi(0) ) )\big) .
$$
A continuous extension of the exponential Euler method on $[0,t_1]$ is given by: for $\theta\in[0,1]$,
$$
U(\theta h_1) = \mathrm{e}^{-\theta h_1A}\phi(0) +  h_1\theta\varphi_1(-\theta h_1A)f\big(0,\phi(0),\phi (-\tau(0,\phi(0) ) )\big) .
$$
For $t\leq 0$, we set $U(t)=\phi(t)$. Once the approximations $u_n\approx u(t_{n})$ and $U(t)\approx u(t)$ in $[0,t_n]$ are obtained, we consider the following local problem
\begin{equation}\label{Eqn:euler-wn+1}
\left\{\begin{aligned}
&w_{n+1}^\prime(t)+Aw_{n+1}(t)=g(t,w_{n+1}(t)),\quad  t\in[t_n,t_{n+1}],\\
&w_{n+1}(0)=u_n,\end{aligned}\right.
\end{equation}
where $g(t,w_{n+1}(t))=f\big(t,w_{n+1}(t),\psi\big(t-\tau(t,w_{n+1}(t)) \big) \big)$ with $\psi$ defined by
\begin{equation*}
\psi(t)=\left\{\begin{aligned} 
&U(t), &&t\in (-\infty,t_n], \\
&w_{n+1}(t), &&t\in[t_n,t_{n+1}]. \end{aligned}\right.
\end{equation*}
Applying the exponential Euler method to \eqref{Eqn:euler-wn+1} leads to
\begin{equation}\label{Eqn:pre-euler-un+1}
u_{n+1} = \mathrm{e}^{-h_{n+1}A}u_n +   h_{n+1}\varphi_1(-  h_{n+1}A)f\big(t_n,u_n,U \big(t_n-\tau(t_n,u_n  ) \big)\big) ,
\end{equation}
where the continuous extension $U(t)$ is already given on $[0,t_n]$. On $[t_n,t_{n+1}]$ it is defined as follows: for $\theta\in[0,1]$,
\begin{equation}\label{Eqn:continuous-euler-Un+1}
U(t_n+\theta h_{n+1}) = \mathrm{e}^{-\theta h_{n+1}A}u_n +    h_{n+1}\theta \varphi_1(-  \theta h_{n+1}A)f\big(t_n,u_n,U \big(t_n-\tau(t_n,u_n  ) \big)\big) .
\end{equation}
The continuous extension satisfies the relations $u_{n}=U(t_{n})$ and $u_{n+1}=U(t_{n+1})$.  

\begin{theorem}
Under the assumptions of Theorem \ref{Thm:well-posedness}, consider for the 
numerical solution of the initial value problem \eqref{Eqn:problem} the exponential Euler method \eqref{Eqn:pre-euler-un+1}-\eqref{Eqn:continuous-euler-Un+1}. For suﬀiciently small $h=\max_{1\leq j \leq N} h_j$, the error bound 
$$
\|u_n-u(t_n)\|\leq Ch 
$$
holds uniformly on $0\leq t \leq T$. The constant $C$ depends on $T$, but is independent of the step size sequence.
\end{theorem}
\begin{proof}
The exact solution of the initial value problem \eqref{Eqn:problem} in $[t_n,t_{n+1}]$ is represented as: for $\theta \in [0,1]$,
$$
\begin{aligned}
u(t_n&+\theta h_{n+1}) =\mathrm{e}^{-\theta h_{n+1}A}u(t_n)\\
& +\int_0^{\theta h_{n+1}} \mathrm{e}^{-(\theta h_{n+1} -\sigma)A} f\big(t_n+\sigma,u(t_n+\sigma),u\big(t_n+\sigma-\tau(t_n+\sigma,u(t_n+\sigma))\big) \big) \,\mathrm{d} \sigma.
\end{aligned}
$$
Denote $e(t)=U(t)-u(t)$. Subtracting the above equality form \eqref{Eqn:continuous-euler-Un+1} gives
\begin{equation}\label{Eqn:e(t)-n+1}
\begin{aligned}
e(t_n&+\theta h_{n+1})  = \mathrm{e}^{-\theta h_{n+1}A}e(t_n)  + R_{n+1}(\theta)\\
&  +  h_{n+1}\theta\varphi_1(-\theta h_{n+1}A)\Big( f\big(t_n,u_n,U\big(t_n-\tau(t_n,u_n) \big) \big) - f\big(t_n,u(t_n),u\big(t_n-\tau(t_n,u(t_n))\big) \big) \Big),
\end{aligned}
\end{equation}
where the local truncation error $R_{n+1}$ is given as
\begin{equation}\label{Eqn:Rn+1}
\begin{aligned}
R_{n+1}(\theta)&=\int_0^{\theta h_{n+1}} \mathrm{e}^{-(\theta h_{n+1}-\sigma)A} \Big( f\big(t_n,u(t_n),u\big(t_n-\tau(t_n,u(t_n))\big) \big)   \\
&\qquad\qquad\qquad\qquad\qquad\qquad- f\big(t_n+\sigma,u(t_n+\sigma),u\big(t_n+\sigma-u(t_n+\sigma)\big) \big) \Big)\,\mathrm{d}\sigma.
\end{aligned}
\end{equation}
Using the boundedness of the semigroup and the Lipschitz condition of $f$, $\tau$ and $u$, one has
$$
\begin{aligned}
\max_{\theta\in[0,1]} \|R_{n+1}\|& \leq C L_f \int_0^{h_{n+1}} \Big(\sigma+\|u(t_n)-u(t_n+\sigma)\| \\
& \quad\qquad\qquad\qquad+\big\|u\big(t_n-\tau(t_n,u(t_n) )\big)-u\big(t_n+\sigma-\tau(t_n+\sigma,u(t_n+\sigma))\big) \big\|\Big)\,\mathrm{d}\sigma \\
&\leq CL_f \int_0^{h_{n+1}} \Big(\sigma+L_u \sigma + L_u \big(\sigma+ L_\tau \sigma+L_\tau L_u  \sigma\big)\Big)\,\mathrm{d}\sigma \leq Ch_{n+1}^2 ,
\end{aligned}
$$
where $L_u= \|u^\prime\|_{L^\infty ( (-\infty,T];X)}$. From \eqref{Eqn:e(t)-n+1} and noting that
\begin{align}
&\big\|U\big(t_{n}-\tau(t_{n},u_{n})\big) - u\big(t_{n}-\tau(t_{n},u(t_{n})) \big) \big\| \notag \\
&\qquad \leq \big\|U\big(t_{n}-\tau(t_{n},u_{n})\big) - u\big(t_{n}-\tau(t_{n},u_{n}) \big) \big\|   + \big\|u\big(t_{n}-\tau(t_{n},u_{n})\big) - u\big(t_{n}-\tau(t_{n},u(t_{n})) \big) \big\| \notag \\
&\qquad \leq \max_{t\leq t_{n}} \|e(t)\| +L_uL_\tau\|e(t_{n})\| ,\label{Eqn:euler-complex-lip}
\end{align}
we obtain 
$$
\begin{aligned}
 \max_{t_n \leq t \leq t_{n+1}} \|e(t)\| 
&  \leq C  \|e(t_n)\|+ Ch_{n+1}^2  + CL_fh_{n+1} \max_{t\leq t_n}\|e(t)\|+CL_f L_u L_\tau h_{n+1} \|e(t_n)\| \\
& \leq C\max_{k=1,\ldots,n} \|e(t_k)\|+ Ch^2+C L_fh \max_{t\leq t_{n+1}}\|e(t)\|.
\end{aligned}
$$
Therefore, for suﬀiciently small $h$, we have
\begin{equation}\label{Eqn:max-etn+1}
\max_{t \leq t_{n+1}} \|e(t)\| \leq C\max_{k=1,\ldots,n} \|e(t_k)\|+ Ch^2 .
\end{equation}

Solving the error recursion \eqref{Eqn:e(t)-n+1} with $\theta=1$ gives
$$
\begin{aligned}
e(t_n) = \sum_{j=1}^{n} \mathrm{e}^{-(t_n-t_j)A}& h_j\varphi_1(-h_jA)\Big( f\big(t_{j-1},u_{j-1},U\big(t_{j-1}-\tau(t_{j-1},u_{j-1})\big) \big)  \\
&   - f\big(t_{j-1},u(t_{j-1}),u\big(t_{j-1}-\tau(t_{j-1},u(t_{j-1})) \big) \big) \Big)  +  \sum_{j=1}^{n} \mathrm{e}^{-(t_n-t_j)A}R_{j}(1),
\end{aligned}
$$
which implies 
$$
\|e(t_n)\|\leq C \sum_{j=1}^n h_j \Big(\|e(t_{j-1})\|+ \big\|U\big(t_{j-1}-\tau(t_{j-1},u_{j-1})\big) - u\big(t_{j-1}-\tau(t_{j-1},u(t_{j-1})) \big) \big\| \Big) +Ch.
$$
Combining the above inequality with \eqref{Eqn:euler-complex-lip} and \eqref{Eqn:max-etn+1}, yields
$$
\|e(t_n)\|\leq C\sum_{j=1}^n h_j   \max_{k=1,\ldots,j-1} \|e(t_k)\|   +Ch.
$$
This further implies that
$$
\max_{k=1,\ldots,n}\|e(t_k)\|\leq C\sum_{j=1}^n h_j   \max_{k=1,\ldots,j-1} \|e(t_k)\|   +Ch.
$$
Applying Gronwall's inequality to the above inequality completes the proof.
\end{proof}

\section{Second order method}\label{Sec:second}
In this section, we construct a second order exponential ERK method for the initial value problem \eqref{Eqn:problem}. Moreover, we establish the well-posedness and convergence of the numerical method. 

As before, the approach consists in solving the local problems step by step and by employing a continuous extension of the numerical solution. For $t\leq 0$, we set $U(t)=\phi(t)$. Once the approximations $u_n\approx u(t_{n})$ and $U(t)\approx u(t)$ in $[0,t_n]$ are obtained, we consider the following local problem
\begin{equation}\label{Eqn:secondERK-wn+1}
\left\{\begin{aligned}
&w_{n+1}^\prime(t)+Aw_{n+1}(t)=g(t,w_{n+1}(t)),\quad  t\in[t_n,t_{n+1}],\\
&w_{n+1}(t_n)=u_n,\end{aligned}\right.
\end{equation}
where $g(t,w_{n+1}(t))=f\big(t,w_{n+1}(t),\psi\big(t-\tau(t,w_{n+1}(t)) \big) \big)$ with $\psi$ defined by
\begin{equation}\label{Eqn:secondERK-xn+1}
\psi(t)=\left\{\begin{aligned} 
&U(t), &&t\in (-\infty,t_n], \\
&w_{n+1}(t), &&t\in[t_n,t_{n+1}]. \end{aligned}\right.
\end{equation}
Applying the second order explicit ERK method \cite[Equation (5.3)]{HochbruckOstermann05:1069} with $c_2\neq 0$ yields 
\begin{equation}\label{Eqn:pre-secondERK-un+1}
\begin{aligned}
\widetilde{u}_{n+1} & = \mathrm{e}^{-h_{n+1}A}u_n +   h_{n+1}\big(\varphi_1(-  h_{n+1} A)-\tfrac{1}{c_2}\varphi_2(-h_{n+1}A) \big) f\big(t_n,u_n,U \big(t_n-\tau(t_n,u_n  ) \big)\big)  \\
&\quad +h_{n+1}\tfrac{1}{c_2} \varphi_2(-h_{n+1}A)f\big(t_{n2},U_{n2},\psi \big(t_{n2}-\tau(t_{n2},U_{n2}  ) \big)\big), \\
U_{n2} & = \mathrm{e}^{-c_2 h_{n+1}A}u_n + c_2  h_{n+1}\varphi_1(-  c_2 h_{n+1} A) f\big(t_n,u_n,U\big(t_n-\tau(t_n,u_n  ) \big)\big),
\end{aligned}
\end{equation}
where $t_{n2}=t_n+c_2h_{n+1}$.

If $\tau$ is bounded from below by a constant $\tau_0>0$, the step size can be choosen as $h_{n+1}\leq \tau_0$, so that the initial value problem \eqref{Eqn:secondERK-wn+1} becomes a problem without delay. As a result, the scheme \eqref{Eqn:pre-secondERK-un+1} with \eqref{Eqn:secondERK-xn+1} is explicit.  However, if $\tau$ can be arbitrary small, then $t_{n2}-\tau(t_{n2},U_{n2})$ may belong to $(t_n,t_{n+1}]$. This phenomenon is referred to as \emph{overlapping}. Since $w_{n+1}$ is unknown, the scheme \eqref{Eqn:pre-secondERK-un+1} with \eqref{Eqn:secondERK-xn+1} is not practical. To address this problem, we construct a continuous numerical solution to approximate $w_{n+1}$. Our starting point is the exact solution of \eqref{Eqn:secondERK-wn+1},  given by
$$
w_{n+1}(t_n+\theta h_{n+1})=\mathrm{e}^{-\theta h_{n+1}A}u_{n}+\int_0^{\theta h_{n+1}} \mathrm{e}^{-(\theta h_{n+1}-\sigma)A}g(t_n+\sigma,w_{n+1}(t_n+\sigma))\,\mathrm{d}\sigma.
$$
The continuous extension $U(t)$ in $[t_n,t_{n+1}]$ is constructed by replacing the term $g(t_n+\sigma,w_{n+1}(t_n+\sigma))$ by the interpolation based on $g(t_n,u_n)$ and $g(t_{n2},U_{n2})$ and replacing $\psi(t)$ by $U(t)$. Consequently, we arrive at the scheme
\begin{equation}\label{Eqn:secondERK-un+1}
\begin{aligned}
u_{n+1} & = \mathrm{e}^{-h_{n+1}A}u_n +   h_{n+1}\big(\varphi_1(-  h_{n+1} A)-\tfrac{1}{c_2}\varphi_2(-h_{n+1}A) \big) f\big(t_n,u_n,U \big(t_n-\tau(t_n,u_n  ) \big)\big)  \\
&\quad +h_{n+1}\tfrac{1}{c_2} \varphi_2(-h_{n+1}A)f\big(t_{n2},U_{n2},U\big(t_{n2}-\tau(t_{n2},U_{n2}  ) \big)\big), \\
U_{n2} & = \mathrm{e}^{-c_2 h_{n+1}A}u_n +   c_2h_{n+1}\varphi_1(-  c_2 h_{n+1} A) f\big(t_n,u_n,U\big(t_n-\tau(t_n,u_n  ) \big)\big),
\end{aligned}
\end{equation}
where $U(t_n+\theta h_{n+1})$, $0\leq \theta \leq 1$ is obtained by
\begin{equation}\label{Eqn:continuous-secondERK-Un+1}
\begin{aligned}
U(t_n+\theta h_{n+1}) & = \mathrm{e}^{-\theta h_{n+1}A}u_n  \\
&\quad +   h_{n+1}\big(\theta \varphi_1(- \theta h_{n+1} A)-\tfrac{1}{c_2}\theta^2\varphi_2(-\theta h_{n+1}A) \big) f\big(t_n,u_n,U \big(t_n-\tau(t_n,u_n  ) \big)\big)  \\
&\quad + h_{n+1}\tfrac{1}{c_2} \theta^2\varphi_2(-\theta h_{n+1}A)f\big(t_{n2},U_{n2},U \big(t_{n2}-\tau(t_{n2},U_{n2}  ) \big)\big).
\end{aligned}
\end{equation}
The continuous extension satisfies $U(t_n)=u_n$ and $U(t_{n+1})=u_{n+1}$, while $U_{n2}\not\equiv U(t_{n2})$. The scheme \eqref{Eqn:secondERK-un+1}-\eqref{Eqn:continuous-secondERK-Un+1} is implicit when overlapping occurs, even if the underlying method is explicit for problem without delay. The following theorem guarantees the unique solvability of the scheme.

\begin{theorem}\label{Thm:well-posedness-second-analytic}
Under the assumptions of Theorem \ref{Thm:well-posedness}, the scheme \eqref{Eqn:secondERK-un+1}-\eqref{Eqn:continuous-secondERK-Un+1} admits a unique solution for sufficiently small step size $h_{n+1}$.
\end{theorem}
\begin{proof}
For $y\in Y =\{v\in C([t_n,t_{n+1}];X) : v(0)=u_n\}$, we define $\widehat{y}:(-\infty,t_{n+1}]\to X$ by
$$
\widehat{y}(t)=\left\{
\begin{aligned}
&U(t), &&t\in (-\infty,t_n],\\
&y(t), &&t\in [t_n,t_{n+1}]. \end{aligned}\right.
$$
We introduce a map $G:Y \to Y$ by: for $t=t_n+\theta h_{n+1}$ with $\theta\in [0,1]$,
\begin{equation*}
\begin{aligned}
G(y)(t)&= \mathrm{e}^{-\theta h_{n+1}A}u_n  \\
&\quad +   h_{n+1}\big(\theta\varphi_1(- \theta h_{n+1} A)-\tfrac{1}{c_2}\theta^2\varphi_2(-\theta h_{n+1}A) \big) f\big(t_n,u_n,U \big(t_n-\tau(t_n,u_n  ) \big)\big)  \\
&\quad + h_{n+1}\tfrac{1}{c_2} \theta^2 \varphi_2(-\theta h_{n+1}A)f\big(t_{n2},U_{n2},\widehat{y} \big(t_{n2}-\tau(t_{n2},U_{n2}  ) \big)\big).
\end{aligned}
\end{equation*}
Using $\|\varphi_2(-tA)\| \leq C_{\mathrm{s}}$ for $t\geq 0$, we obtain that, for $y_1,y_2\in C([t_n,t_{n+1}];X)$,
$$
\begin{aligned}
\|G(y_1)-G(y_2)\|_{C([t_n,t_{n+1}];X)} \leq C_{\mathrm{s}}L_f h_{n+1}\tfrac{1}{c_2} \|y_1-y_2\|_{C([t_n,t_{n+1}];X)},
\end{aligned}
$$
which implies that $G$ is a contraction on $C([t_n,t_{n+1}];X)$ for $h_{n+1}< c_2(C_{\mathrm{s}}L_f)^{-1}$. 
It follows from the Banach fixed point theorem that the equation $G(y)=y$ has a unique solution, which completes the proof.
\end{proof}

For the error analysis, we introduce the local problem 
\begin{equation}\label{Eqn:secondERK-zn+1}
\left\{\begin{aligned}
&z_{n+1}^\prime(t)+Az_{n+1}(t)=f\big(t,z_{n+1}(t),u\big(t-\tau(t,z_{n+1}(t)) \big) \big),\quad  t\in[t_n,t_{n+1}],\\
&z_{n+1}(t_{n})=u(t_n),\end{aligned}\right.
\end{equation}
whose solution obviously is $z(t)=u(t)$. Consider for its numerical solution 
\begin{equation}\label{Eqn:secondERK-hatun+1}
\begin{aligned}
\widehat{u}_{n+1} & = \mathrm{e}^{-h_{n+1}A}u(t_n) +   h_{n+1}\big(\varphi_1(-  h_{n+1} A)-\tfrac{1}{c_2}\varphi_2(-h_{n+1}A) \big) f\big(t_n,u(t_n),u \big(t_n-\tau(t_n,u(t_n)  ) \big)\big)  \\
&\quad +h_{n+1}\tfrac{1}{c_2} \varphi_2(-h_{n+1}A)f\big(t_{n2},\widehat{U}_{n2},u \big(t_{n2}-\tau(t_{n2},\widehat{U}_{n2}  ) \big)\big), \\
\widehat{U}_{n2} & = \mathrm{e}^{-c_2 h_{n+1}A}u(t_n) +   c_2h_{n+1}\varphi_1(-  c_2 h_{n+1} A) f\big(t_n,u(t_n),u\big(t_n-\tau(t_n,u(t_n)  ) \big)\big)
\end{aligned}
\end{equation}
and the corresponding continuous extension
\begin{equation}\label{Eqn:continuous-secondERK-hatUn+1}
\begin{aligned}
\widehat{U}(t_n+\theta h_{n+1}) & = \mathrm{e}^{-\theta h_{n+1}A}u(t_n)  \\
&\quad +   h_{n+1}\big(\theta\varphi_1(- \theta h_{n+1} A)-\tfrac{1}{c_2}\theta^2\varphi_2(-\theta h_{n+1}A) \big) f\big(t_n,u(t_n),u \big(t_n-\tau(t_n,u(t_n)  ) \big)\big)  \\
&\quad + h_{n+1}\tfrac{1}{c_2} \theta^2\varphi_2(-\theta h_{n+1}A)f\big(t_{n2},\widehat{U}_{n2},u \big(t_{n2}-\tau(t_{n2},\widehat{U}_{n2}  ) \big)\big).
\end{aligned}
\end{equation}
The local error estimate is given in the next lemma.

\begin{lemma}\label{Lem:second-order}
Under the Assumptions \ref{Ass:sectorial}-\ref{Ass:lip}, if the function
$$
g(t)=f\big(t,u(t),u\big(t-\tau(t,u(t))\big)\big)
$$
is of class $C^{1,1}$ on $[t_n ,t_{n+1}]$, then the following error bounds
$$
\begin{aligned}
\|\widehat{U}_{n2}-u(t_{n2})\|&\leq Ch_{n+1}^2, \\
\max_{t_n\leq t \leq t_{n+1}} \|\widehat{U}(t)-u(t)\|&\leq Ch_{n+1}^3,
\end{aligned}
$$
hold. The constant $C$ is independent of $h_{n+1}$.
\end{lemma}
\begin{proof}
Expanding $g$ into a Taylor series with remainder in integral form, the solution $u$ on $[t_n,t_{n+1}]$ can be written as
\begin{align}
u(t_n+\theta h_{n+1})&=\mathrm{e}^{-\theta h_{n+1}A}u(t_n)+
\int_0^{\theta h_{n+1}}\mathrm{e}^{-(\theta h_{n+1}-\sigma)A}g(t_n+\sigma)\,\mathrm{d}\sigma \notag \\
&=\mathrm{e}^{-\theta h_{n+1}A}u(t_n)+ \theta h_{n+1}\varphi_1(-\theta h_{n+1}A)g(t_n) + (\theta h_{n+1} )^2 \varphi_2(-\theta h_{n+1}A)g^\prime(t_n)  \notag \\
&\quad+ \int_0^{\theta h_{n+1}}\mathrm{e}^{-(\theta h_{n+1}-\sigma)A}\int_0^\sigma (\sigma - \xi)g^{(2)}(t_n+\xi)\,\mathrm{d}\xi\,\mathrm{d}\sigma.\label{Eqn:taylor-exactun+1}
\end{align}
Since $g \in C^{1,1}([t_n, t_{n+1}]; X)$, its second derivative $g^{(2)}$ exists almost everywhere on $(t_n, t_{n+1})$ satisfying $g^{(2)} \in L^\infty(t_n, t_{n+1}; X)$. On the other hand, plugging the solution $u$ into \eqref{Eqn:continuous-secondERK-hatUn+1} (with $\widehat{U}$ replaced by $u$ and $\widehat{U}_{n2}$ replaced by $u(t_{n2}))$ gives  
\begin{equation*} 
\begin{aligned}
u(t_n+\theta h_{n+1}) & = \mathrm{e}^{-\theta h_{n+1}A}u(t_n)   +   h_{n+1}\big(\theta\varphi_1(- \theta h_{n+1} A)-\tfrac{1}{c_2}\theta^2\varphi_2(-\theta h_{n+1}A) \big) g(t_n) \\
&\quad + h_{n+1}\tfrac{1}{c_2} \theta^2 \varphi_2(-\theta h_{n+1}A)g(t_{n2}) + \Delta_{n+1}(\theta),
\end{aligned}
\end{equation*}
with defect $\Delta_{n+1}(\theta)$. Now expanding $g$ into a Taylor series with remainder in integral form gives
\begin{equation}\label{Eqn:taylor-continuous-secondERK-exactun+1}
\begin{aligned}
u(t_n+\theta h_{n+1}) & = \mathrm{e}^{-\theta h_{n+1}A}u(t_n)   +  \theta h_{n+1} \varphi_1(- \theta h_{n+1} A) g(t_n)  + (\theta h_{n+1})^2 \varphi_2(-\theta h_{n+1}A) g^\prime(t_{n})  \\
&\quad+ h_{n+1}\tfrac{1}{c_2} \theta^2 \varphi_2(-\theta h_{n+1}A) \int_0^{c_2h_{n+1}} (c_2h_{n+1} - \sigma)g^{(2)}(t_n+\sigma)\,\mathrm{d}\sigma + \Delta_{n+1}(\theta).
\end{aligned}
\end{equation}
Subtracting \eqref{Eqn:taylor-exactun+1} from \eqref{Eqn:taylor-continuous-secondERK-exactun+1} gives the following explicit representation of the defect,
$$
\begin{aligned}
\Delta_{n+1}(\theta) & =  \int_0^{\theta h_{n+1}}\mathrm{e}^{-(\theta h_{n+1}-\sigma)A}\int_0^\sigma (\sigma - \xi)g^{(2)}(t_n+\xi)\,\mathrm{d}\xi\,\mathrm{d}\sigma\\
&\quad- h_{n+1}\tfrac{1}{c_2} \theta^2 \varphi_2(-\theta h_{n+1}A) \int_0^{c_2h_{n+1}} (c_2h_{n+1} - \sigma)g^{(2)}(t_n+\sigma)\,\mathrm{d}\sigma ,
\end{aligned}
$$
which implies
$$
\max_{\theta\in [0,1]} \|\Delta_{n+1}(\theta)\| \leq Ch_{n+1}^3 \|g^{(2)}\|_{L^\infty([t_n,t_{n+1}];X)}.
$$
Finally, noting that
$$
\begin{aligned}
&\widehat{U}(t_n+\theta h_{n+1}) - u(t_n+\theta h_{n+1}) \\
&\qquad=
  h_{n+1}\tfrac{1}{c_2} \theta^2\varphi_2(-\theta h_{n+1}A) \Big(f\big(t_{n2},\widehat{U}_{n2},u \big(t_{n2}-\tau(t_{n2},\widehat{U}_{n2}  ) \big)\big) - g(t_{n2}) \Big) - \Delta_{n+1}(\theta)
\end{aligned}
$$
and using \eqref{Eqn:secondERK-hatun+1} and \eqref{Eqn:Rn+1}, which shows
$$
\|\widehat{U}_{n2}-u(t_{n2})\|=\|R_{n+1}(c_2)\| \leq Ch_{n+1}^2,
$$
we obtain  
$$
\max_{\theta\in[0,1]}\|\widehat{U}(t_n+\theta h_{n+1}) - u(t_n+\theta h_{n+1})\| \leq Ch_{n+1}^3.
$$
Thus, the proof is completed.
\end{proof}

Let $\widehat{e}(t)=U(t)-\widehat{U}(t)$ and $\widetilde{e}(t)=\widehat{U}(t)-u(t)$. Then, we have
$$
e(t)=U(t)-u(t)=\widehat{e}(t)+\widetilde{e}(t).
$$ 
From \eqref{Eqn:secondERK-un+1} and
$$
u(t_{n2})  = \mathrm{e}^{-c_2 h_{n+1}A}u(t_n)   +   c_2 h_{n+1} \varphi_1(- c_2 h_{n+1} A)  g(t_n) + Ch_{n+1}^2,
$$
we obtain
\begin{equation}\label{Eqn:Un2-u(tn2)}
\|U_{n2}-u(t_{n2})\|  \leq C\|e(t_n)\| + Ch_{n+1}  \max_{t\leq t_n}\|e(t)\|+Ch_{n+1}^2.
\end{equation}
Now, we are ready to prove the convergence of the scheme \eqref{Eqn:secondERK-un+1}-\eqref{Eqn:continuous-secondERK-Un+1} under the assumption that
\begin{equation}\label{Eqn:g-C11}
g(t)=f\big(t,u(t),u\big(t-\tau(t,u(t))\big)\big) \in C^{1,1}([t_j,t_{j+1}];X),\quad j=0,\ldots,N-1.
\end{equation}

\begin{theorem}\label{Thm:second-order}
Under the Assumptions \ref{Ass:sectorial}-\ref{Ass:lip}, let $g$ satisfy the condition \eqref{Eqn:g-C11}. Consider for the numerical solution of the initial value problem \eqref{Eqn:problem} the second order ERK method \eqref{Eqn:secondERK-un+1}-\eqref{Eqn:continuous-secondERK-Un+1}. For suﬀiciently small $h$, the error bound
$$
\|u_n-u(t_n)\|\leq Ch^2 
$$
holds uniformly on $0\leq t \leq T$. The constant $C$ depends on $T$, but is independent of the step size sequence.
\end{theorem}
\begin{proof}
Subtracting \eqref{Eqn:continuous-secondERK-hatUn+1} from \eqref{Eqn:continuous-secondERK-Un+1} yields
\begin{equation}\label{Eqn:second-recursion}
\begin{aligned}
e&(t_n+\theta h_{n+1})  = \mathrm{e}^{-\theta h_{n+1}A}e(t_n) +\widetilde{e}(t_n+\theta h_{n+1})  \\
&\quad +   h_{n+1}b_1(\theta;-h_{n+1}A)  \Big( f\big(t_n,u_n,U \big(t_n-\tau(t_n,u_n  ) \big)\big) - f\big(t_n,u(t_n),u \big(t_n-\tau(t_n,u(t_n)  ) \big)\big) \Big) \\
&\quad + h_{n+1}b_2(\theta;-h_{n+1}A)\Big( f\big(t_{n2},U_{n2},U \big(t_{n2}-\tau(t_{n2},U_{n2}  ) \big)\big) - f\big(t_{n2},\widehat{U}_{n2},u \big(t_{n2}-\tau(t_{n2},\widehat{U}_{n2}  ) \big)\big) \Big),
\end{aligned}
\end{equation}
where the weights $b_i(\theta;-tA)$ are defined by
\begin{equation}\label{Eqn:b1b2}
b_1(\theta;-h_{n+1}A)=\theta\varphi_1(- \theta h_{n+1} A)-\tfrac{1}{c_2}\theta^2\varphi_2(-\theta h_{n+1}A) , \quad
b_2(\theta;-h_{n+1}A)=\tfrac{1}{c_2} \theta^2\varphi_2(-\theta h_{n+1}A).
\end{equation} 
Noting that
\begin{align}
&\big\|U \big(t_n-\tau(t_n ,u_n  ) \big)-u \big(t_n-\tau(t_n,u(t_n)  ) \big)\big\| \notag\\
& \qquad\qquad\qquad \leq \|U \big(t_n-\tau(t_n,u_n  ) \big)\big)-u \big(t_n-\tau(t_n,u_n  ) \big)\big\| \notag\\
& \qquad\qquad\qquad\quad+\|u \big(t_n-\tau(t_n,u_n  ) \big)-u \big(t_n-\tau(t_n,u(t_n)  ) \big)\big\| \notag\\
& \qquad\qquad\qquad \leq  \max_{t\leq t_n}\|e(t)\| + C\|e(t_n)\|,\label{Eqn:second-complicanted-lip-1} \\
& \big\|U \big(t_{n2}-\tau(t_{n2}  ,U_{n2}  ) \big) -u \big(t_{n2}-\tau(t_{n2},\widehat{U}_{n2}  ) \big)\big\| \notag\\
&  \qquad\qquad\qquad\leq  \big\|U \big(t_{n2}-\tau(t_{n2},U_{n2}  ) \big) -u \big(t_{n2}-\tau(t_{n2},U_{n2}  ) \big)\big\|  \notag\\
& \qquad\qquad\qquad\quad   +\big\|u \big(t_{n2}-\tau(t_{n2},U_{n2}  ) \big) -u \big(t_{n2}-\tau(t_{n2},u(t_{n2})  ) \big)\big\| \notag\\
& \qquad\qquad\qquad\quad+\big\|u \big(t_{n2}-\tau(t_{n2},u(t_{n2})  ) \big) -u \big(t_{n2}-\tau(t_{n2},\widehat{U}_{n2}  ) \big)\big\| \notag\\
& \qquad\qquad\qquad \leq \max_{t\leq t_n} \|e(t)\|+C\|U_{n2}-u(t_{n2})\|+C\|u(t_{n2})-\widehat{U}_{n2}\|
\label{Eqn:second-complicanted-lip-2}
\end{align}
and using Lemma \ref{Lem:second-order} and \eqref{Eqn:Un2-u(tn2)}, we obtain from \eqref{Eqn:second-recursion}
$$
\begin{aligned}
&\max_{t_n\leq t \leq t_{n+1}} \|e(t)\| \\
&\qquad\leq C\|e(t_n)\| +Ch_{n+1}^3 +Ch_{n+1}\Big(\|e(t_n)\| +\big\|U \big(t_n-\tau(t_n,u_n  ) \big)\big)-u \big(t_n-\tau(t_n,u(t_n)  ) \big)\big\| \Big) \\
&\qquad\quad+ Ch_{n+1}\Big( \|U_{n2}-\widehat{U}_{n2}\| +\big\|U \big(t_{n2}-\tau(t_{n2},U_{n2}  ) \big) -u \big(t_{n2}-\tau(t_{n2},\widehat{U}_{n2}  ) \big)\big\|\Big) \\
&\qquad\leq C\max_{k=1,\ldots,n}\|e(t_k)\| +Ch^3 +Ch \max_{t\leq t_{n+1}}\|e(t)\|      .
\end{aligned}
$$
Therefore, for suﬀiciently small $h$, we have
\begin{equation}\label{Eqn:second-maxtn+1-e(t)}
\max_{  t \leq t_{n+1}} \|e(t)\| \leq C\max_{k=1,\ldots,n}\|e(t_k)\| +Ch^3 .
\end{equation}

Solving the recursion \eqref{Eqn:second-recursion} with $\theta=1$ gives
$$
\begin{aligned}
e(t_n) & =    \sum_{j=1}^n\mathrm{e}^{-(t_n-t_j)A} h_{j} \bigg[ b_1(1;-h_{j}A)  \Big( f\big(t_{j-1},u_{j-1},U \big(t_{j-1}-\tau(t_{j-1},u_{j-1}  ) \big)\big)\\
&\qquad\qquad\qquad\qquad\quad     - f\big(t_{j-1},u(t_{j-1}),u \big(t_{j-1}-\tau(t_{j-1},u(t_{j-1})  ) \big)\big) \Big) \\
&\qquad\qquad\qquad\qquad\quad  +  b_2(1;-h_{j}A)\Big( f\big(t_{j-1,2},U_{j-1,2},U \big(t_{j-1,2}-\tau(t_{j-1,2},U_{j-1,2}  ) \big)\big) \\
&\qquad\qquad\qquad\qquad\quad   - f\big(t_{j-1,2},\widehat{U}_{j-1,2},u \big(t_{j-1,2}-\tau(t_{j-1,2},\widehat{U}_{j-1,2}  ) \big)\big) \Big) \bigg]+ \sum_{j=1}^n\mathrm{e}^{-(t_n-t_j)A}\widetilde{e}(t_j) ,
\end{aligned}
$$
which, together with Lemma \ref{Lem:second-order} and \eqref{Eqn:second-complicanted-lip-1}, \eqref{Eqn:second-complicanted-lip-2}, \eqref{Eqn:second-maxtn+1-e(t)}, implies
$$
\begin{aligned}
\|e(t_n)\| &\leq C\sum_{j=1}^n h_{j} \max_{k=1,\ldots,j-1} \|e(t_k)\|  +Ch^3 +\Bigg\|\sum_{j=1}^n\mathrm{e}^{-(t_n-t_j)A}\widetilde{e}(t_j)\Bigg\|\\
&\leq C\sum_{j=1}^n h_{j} \max_{k=1,\ldots,j-1} \|e(t_k)\|  +Ch^2.
\end{aligned} 
$$
This further implies that
$$
\max_{k=1,\ldots,n}\|e(t_k)\|\leq C\sum_{j=1}^n h_j   \max_{k=1,\ldots,j-1} \|e(t_k)\|   +Ch^2.
$$
Applying Gronwall's inequality to the above inequality completes the proof.
\end{proof}
\begin{remark}\label{Rem:second}\rm
Let $f$ and $\tau$ be of class $C^{1,1}$ on their respective domains, and $\phi \in C^{1,1}((-\infty, T]; X)$ with $\phi(0) \in D(A)$. If $-A\phi(0) + f\big(0, \phi(0), \phi\big(0 - \tau(0, \phi(0))\big)\big) = \phi'(0) \in D(A)$, then the solution $u \in C^{1,1}((-\infty, T]; X)$ of the initial value problem \eqref{Eqn:problem} is guaranteed by Theorem \ref{Thm:well-posedness}. It follows that $g \in C^{1,1}((-\infty, T]; X)$ and thereby the condition \eqref{Eqn:g-C11} holds for arbitrary meshes.

On the other hand, if $-A\phi(0) + f\big(0, \phi(0), \phi(0 - \tau(0, \phi(0)))\big) \neq \phi'(0)$, assuming that the solution satisfies $u\in C^1([0,T];D(A))$, it follows that the solution $u$ belongs to $C^{1,1}( [0,T];X)$ but possesses a discontinuity in its derivative at $t=0$; that is, $u^{\prime}(0^+)\neq\phi^{\prime}(0^-)$. As a result, the function $g\not\in C^{1,1}( [t_j,t_{j+1}];X)$ if $0 \in  \{ t-\tau(t,u(t)): t \in (t_j,t_{n+1})\}$. However, second order convergence typically holds even for arbitrarily meshes. In particular, we consider the case where the set of discontinuities $M=\{ t \in (0,T) : t -\tau(t,u(t))=0\}$ has only a few elements. Given a mesh $I_h$, we define the set of indices corresponding to mesh intervals that contain at least one discontinuity as $J_M = \{ j : M \cap (t_{j-1},t_{j}) \neq \emptyset,\, j=1,\ldots,N\}$. Then it holds
$$
\Bigg\|\sum_{j=1}^n\mathrm{e}^{-(t_n-t_j)A}\widetilde{e}(t_j)\Bigg\|\leq 
\Bigg\|\sum_{\substack{j=1 \\ j\notin J_M}}^n\mathrm{e}^{-(t_n-t_j)A}\Delta_{j}(1)\Bigg\| 
+ \Bigg\|\sum_{j\in J_M}\mathrm{e}^{-(t_n-t_j)A}R_j(1)\Bigg\| \leq Ch^2.
$$
It follows that the second order convergence result remains valid.

\end{remark}

It is possible to develop higher order ERK methods for \eqref{Eqn:problem}. In particular, arbitrary high order methods can be systematically constructed by using the methods of collocation type \cite{HochbruckOstermann05:323}.

\section{Higher order methods of collocation type}\label{Sec:collocation}
In this section, we extend the ERK methods of collocation type for (non\-vanishing) time-dependent delay, as developed in \cite{HuangOstermann25}, to the initial value problem \eqref{Eqn:problem} with arbitrary state-dependent delay.

For $t\leq 0$, we set $U(t)=\phi(t)$. Once the approximations $u_n\approx u(t_{n})$ and $U(t)\approx u(t)$ in $[0,t_n]$ are obtained, again, we consider the local problem
\begin{equation}\label{Eqn:colERK-wn+1}
\left\{\begin{aligned}
&w_{n+1}^\prime(t)+Aw_{n+1}(t)=g(t,w_{n+1}(t)),\quad  t\in[t_n,t_{n+1}],\\
&w_{n+1}(t_n)=u_n,\end{aligned}\right.
\end{equation}
where $g(t,w_{n+1}(t))=f\big(t,w_{n+1}(t),\psi \big(t-\tau(t,w_{n+1}(t)) \big) \big)$ with $\psi$ defined by
\begin{equation}\label{Eqn:col-local-x(s)}
\psi(t)=\left\{\begin{aligned} 
&U(t), &&t\in (-\infty,t_n], \\
&w_{n+1}(t), &&t\in[t_n,t_{n+1}]. \end{aligned}\right.
\end{equation}
The solution can be represented as
\begin{equation}\label{Eqn:col-local-problem}
w_{n+1}(t_n+\theta h_{n+1})=\mathrm{e}^{-\theta h_{n+1}A}u_n +
\int_0^{\theta h_{n+1}}\mathrm{e}^{-(\theta h_{n+1}-\sigma)A}g(t_n+\sigma,w_{n+1}(t_n+\sigma))\,\mathrm{d}\sigma.
\end{equation}
Applying the ERK method of collocation type \cite[Equation (4)]{HochbruckOstermann05:323} with nonconfluent nodes $c_1,\ldots,c_s\in[0,1]$, yields 
\begin{equation}\label{Eqn:pre-colERK-un+1}
\begin{aligned}
\widetilde{u}_{n+1} & = \mathrm{e}^{-h_{n+1}A}u_n +   h_{n+1}\sum_{i=1}^s b_i(-h_{n+1}A) f\big(t_{ni},\widetilde{U}_{ni},\psi \big(t_{ni}-\tau(t_{ni},\widetilde{U}_{ni}  ) \big)\big),  \\
\widetilde{U}_{ni} & = \mathrm{e}^{-c_i h_{n+1}A}u_n +h_{n+1} \sum_{j=1}^s a_{ij}(-h_{n+1}A) f\big(t_{nj},\widetilde{U}_{ni},\psi \big(t_{nj}-\tau(t_{nj},\widetilde{U}_{nj}  ) \big)\big),\quad 1\leq i \leq s,
\end{aligned}
\end{equation}
where $t_{ni}=t_n+c_i h_{n+1}$ and 
$$
\begin{aligned}
a_{ij}(-h_{n+1} A) & = \frac{1}{h_{n+1}} \int_0^{c_i h_{n+1}} 
\mathrm{e}^{-(c_i h_{n+1} - \sigma) A} \, \ell_j\bigg( \frac{\sigma}{h_{n+1}} \bigg) \, \mathrm{d}\sigma, \\
b_{i}(-h_{n+1} A) & = \frac{1}{h_{n+1}} \int_0^{  h_{n+1}} 
\mathrm{e}^{-(  h_{n+1} - \sigma) A} \, \ell_i\bigg( \frac{\sigma}{h_{n+1}} \bigg) \, \mathrm{d} \sigma \quad\mbox{with}~\ell_i(\rho) = \prod_{m = 1 , m \ne i}^{s} 
\frac{\rho - c_m}{c_i - c_m}.
\end{aligned}
$$
When overlapping occurs, the scheme \eqref{Eqn:pre-colERK-un+1} with \eqref{Eqn:col-local-x(s)} is not practical since $w_{n+1}$ is unkonwn. We therefore modify this scheme and construct the continuous numerical solution $U(t)$ in $[t_n,t_{n+1}]$ by replacing the term $g(t_n+\sigma,w_{n+1}(t_n+\sigma))$ in \eqref{Eqn:col-local-problem} by the interpolation based on $g(t_{ni},U_{ni})$ ($i=1,\ldots,s$) and replacing $\psi(t)$ by $U(t)$. Consequently, we arrive at the scheme
\begin{equation}\label{Eqn:colERK-un+1}
\begin{aligned}
{u}_{n+1} & = \mathrm{e}^{-h_{n+1}A}u_n +   h_{n+1}\sum_{i=1}^s b_i(-h_{n+1}A) f\big(t_{ni},{U}_{ni},U\big(t_{ni}-\tau(t_{ni},{U}_{ni}  ) \big)\big),  \\
{U}_{ni} & = \mathrm{e}^{-c_i h_{n+1}A}u_n +h_{n+1} \sum_{j=1}^s a_{ij}(-h_{n+1}A) f\big(t_{nj},{U}_{ni},U\big(t_{nj}-\tau(t_{nj},{U}_{nj}  ) \big)\big),\quad 1\leq i \leq s,
\end{aligned}
\end{equation}
where $U(t_n+\theta h_{n+1})$, $0\leq \theta \leq 1$, is defined by
\begin{equation}\label{Eqn:continuous-colERK-Un+1}
\begin{aligned}
U(t_n+\theta h_{n+1}) & = \mathrm{e}^{-\theta h_{n+1}A}u_n + h_{n+1}\sum_{i=1}^s b_i(\theta;-h_{n+1}A) f\big(t_{ni},{U}_{ni},U \big(t_{ni}-\tau(t_{ni},{U}_{ni}  ) \big)\big),
\end{aligned}
\end{equation}
and
$$
b_{i}(\theta;-h_{n+1} A)  = \frac{1}{h_{n+1}}  \int_0^{  \theta h_{n+1}} 
\mathrm{e}^{-( \theta h_{n+1} - \sigma) A} \, \ell_j\bigg( \frac{\sigma}{h_{n+1}} \bigg) \, \mathrm{d}\sigma.
$$
Recalling the boundedness of $\mathrm{e}^{tA}$, we have the estimate
\begin{equation}\label{Eqn:boundedness-bitheta}
\| b_i(\theta;- t A) \|_{X\leftarrow X} \leq C,\quad  t\geq0.
\end{equation}
Noting the relations 
$$
b_{i}(1;-h_{n+1} A)=b_{i}( -h_{n+1} A),\quad 
b_{j}(c_i;-h_{n+1} A)=a_{ij}(-h_{n+1} A),\quad 1\leq i,j \leq s,
$$
we obtain $U(t_n)=u_n$, $U(t_{n+1})=u_{n+1}$ and $U(t_{ni})=U_{ni}$ for $i=1,\ldots,s$. The scheme \eqref{Eqn:colERK-un+1}-\eqref{Eqn:continuous-colERK-Un+1} remains implicit regardless of whether overlapping occurs or not (except for the case where $s=1$ and $c_1=0$). The following theorem guarantees the unique solvability of the scheme.

\begin{theorem}\label{Thm:well-posedness-col}
Under the assumptions of Theorem \ref{Thm:well-posedness}, the scheme \eqref{Eqn:colERK-un+1}-\eqref{Eqn:continuous-colERK-Un+1} admits a unique solution for sufficiently small step size $h_{n+1}$.
\end{theorem}
\begin{proof}
We first show that the local problem \eqref{Eqn:colERK-wn+1} with \eqref{Eqn:col-local-x(s)} is well-posed. On each interval $[t_j,t_{j+1}]$ the continuous extension $U(t)$ is represented as
$$
U(t)=\mathrm{e}^{-(t-t_j)A}u_j + \int_0^{t-t_j}\mathrm{e}^{-(t-t_j-\sigma)A} r_{j+1}(\sigma) \,\mathrm{d}\sigma, \quad t_j\leq t \leq t_{j+1},
$$
where $r_{j+1} \in C([t_j,t_{j+1}];X)$ has the form
$$
r_{j+1}(\sigma)=\sum_{i=1}^s \ell_i \bigg( \frac{\sigma}{h_{j+1}}\bigg) f\big(t_{ji},U_{ji},U\big(t_{ji}-\tau(t_{ji},U_{ji})\big)\big) .
$$
This a direct consequence of \eqref{Eqn:continuous-colERK-Un+1}.
Note that $U(t)$ is thus the solution of the initial value problem
$$
\left\{\begin{aligned}
&W^\prime(t)+AW(t)= r_{j+1}(t), && t\in[t_j,t_{j+1}], \\
&W(t)=U(t),&& t\in(-\infty,t_j].
\end{aligned}\right.
$$
An application of \cite[Theorem 4.3.1]{Lunardi95Book} yields $U\in C^1([t_j,t_{j+1}];X)\cap C([t_j,t_{j+1}];D(A))$. As a result, $U(t)$ is Lipschitz continuous for $t\leq t_n$ and $u_n\in D(A)$. By Theorem \ref{Thm:well-posedness}, the problem \eqref{Eqn:colERK-wn+1} with \eqref{Eqn:col-local-x(s)} admits a unique solution $w_{n+1}$ up to time $t_{n+1}$ by choosing $h_{n+1}$ sufficiently small.

For $y\in Y =\{v\in C([t_n,t_{n+1}];X) : v(0)=u_n\}$, we define $\widehat{y}:(-\infty,t_{n+1}]\to X$ by the equation
$$
\widehat{y}(t)=\left\{
\begin{aligned}
&U(t), &&t\in (-\infty,t_n],\\
&y(t), &&t\in [t_n,t_{n+1}]  \end{aligned}\right.
$$
and introduce a map $G:Y \to Y$ by: for $t=t_n+\theta h_{n+1}$ with $\theta\in [0,1]$,
\begin{equation}\label{Eqn:mapping-g}
\begin{aligned}
G(y)(t)&= \mathrm{e}^{-\theta h_{n+1}A}u_n + h_{n+1}\sum_{i=1}^s b_i(\theta;-h_{n+1}A) f\big(t_{ni},y(t_{ni}),\widehat{y} \big(t_{ni}-\tau(t_{ni},y(t_{ni})  ) \big)\big)  .
\end{aligned}
\end{equation}
Let $B = \{y \in Y : \|y- {G(w_{n+1})}\|_{C([t_n,t_{n+1}];X)} \leq 1\}$. The set $B$ is a nonempty, closed, bounded, convex subset of $C([t_n,t_{n+1}];X)$. Using $\|b_i(\theta;- h_{n+1}A)\|_{X\leftarrow X} \leq C_{\mathrm{s}}$ and the fact that
\begin{align}
& \big\|f\big(t_{ni},w_{n+1}(t_{ni}),\widehat{w}_{n+1} \big(t_{ni}-\tau(t_{ni},w_{n+1}(t_{ni})  ) \big)\big)\big\| \notag \\
&\qquad\leq \big\|f\big(t_{ni},w_{n+1}(t_{ni}),\widehat{w}_{n+1} \big(t_{ni}-\tau(t_{ni},w_{n+1}(t_{ni})  ) \big)\big)   - f(0,0,0) \big\| +  \|f(0,0,0) \| \notag \\
&\qquad \leq L_fT + 2L_f\|w_{n+1}\|_{C([t_n,t_{n+1}];X)} + L_f \|U\|_{C( (-\infty,t_n];X)}+  \|f(0,0,0) \|,\notag
\end{align}
we have
\begin{align}
&\big\|G(y) - {G(w_{n+1})}\big\|_{C([t_n,t_{n+1}];X)} \notag \\
&\quad \leq  C_fC_{\mathrm{s}} h_{n+1} \sum_{i=1}^s \Big(\|y(t_{ni})-w(t_{ni})\|+ \big\| \widehat{y} \big(t_{ni}-\tau(t_{ni},y(t_{ni})  ) \big)- \widehat{w}_{n+1}\big(t_{ni}-\tau(t_{ni},w(t_{ni})  ) \big)\big\|\Big) \notag \\
&\quad \leq  C_fC_{\mathrm{s}} h_{n+1}s  \Big(\|y - {G(w_{n+1})}\|_{C([t_n,t_{n+1}];X)}+  \| w_{n+1}  \|_{C([t_n,t_{n+1}];X)} + \|   {G(w_{n+1})} \|_{C([t_n,t_{n+1}];X)}  \Big) \notag \\
&\qquad+   C_fC_{\mathrm{s}} h_{n+1} \sum_{i=1}^s \Big(  \big\| \widehat{y} \big(t_{ni}-\tau(t_{ni},y(t_{ni})  ) \big)- \widehat{G(w_{n+1})}\big(t_{ni}-\tau(t_{ni},y(t_{ni})  ) \big)\big\|   \Big) \notag \\
&\qquad+   C_fC_{\mathrm{s}} h_{n+1} \sum_{i=1}^s \Big(  \big\|\widehat{G(w_{n+1})}\big(t_{ni}-\tau(t_{ni},y(t_{ni})  ) \big)\big\| + \big\| \widehat{w}_{n+1}\big(t_{ni}-\tau(t_{ni},w(t_{ni})  ) \big)\big\|   \Big) \notag \\
&\quad \leq Ch_{n+1}\big( \|U\|_{C( (-\infty,t_n];X)} + \|w_{n+1}\|_{C([t_n,t_{n+1}];X)} +  \|f(0,0,0) \| +1   \big).
\notag
\end{align}
For $h_{n+1}$ sufficiently small, this shows that
$$
\big\|G(y) - {G(w_{n+1})}\big\|_{C([t_n,t_{n+1}];X)} \leq 1.
$$
Thus $G$ maps $B$ to $B$. 

Inspired by \cite{{Fitzgibbon78:1}}, the existence of a solution to the scheme \eqref{Eqn:colERK-un+1}-\eqref{Eqn:continuous-colERK-Un+1} can be established by Schauder's fixed point theorem \cite[p.~179]{Brezis11Book} applied to the map $G$. The continuity of $G$ is straightforward to verify. It remains to show that $G(B)$ is precompact, i.e., its closure is compact in $X$. For this purpose, we first show that for all $t=t_n+\theta h_{n+1} \in [t_n,t_{n+1}]$, the set $\{G(y)(t):y\in B\}$ is precompact. For an arbitrary $\theta\in(0,1]$, we choose $\xi \in (0,\theta)$. For $y\in B$, we define 
$$
\begin{aligned}
G_{\xi}(y)(t) & =\mathrm{e}^{-\theta h_{n+1} A}u_n  +    \sum_{i=1}^s\int_0^{(\theta-\xi) h_{n+1}} \mathrm{e}^{- (\theta  h_{n+1}-\sigma)A}K_y^{(i)}(\sigma) \, \mathrm{d}\sigma  \\
& = \mathrm{e}^{-\theta h_{n+1} A}u_n  +    \sum_{i=1}^s \mathrm{e}^{-\xi h_{n+1}A}\int_0^{(\theta-\xi) h_{n+1}} \mathrm{e}^{- ( (\theta-\xi)  h_{n+1}-\sigma)A}K_y^{(i)}(\sigma) \, \mathrm{d}\sigma .
\end{aligned}
$$
where
$$
K_y^{(i)}(\sigma)= \ell_i \bigg( \frac{\sigma}{h_{n+1}}\bigg) f\big(t_{ni},y(t_{ni}),\widehat{y}\big(t_{ni}-\tau(t_{ni},y(t_{ni}))\big)\big).
 $$
Since $\mathrm{e}^{-tA}$ is compact for $t>0$, the set $\{G_{\xi}(y)(t):y\in B\}$ is precompact in $X$. Noting that $G_{\xi}(y)(t)\to G(y)(t)$ in $X$ as $\xi\to 0$, we obtain that $\{G(y)(t):y\in B\}$ is totally bounded and thereby precompact. We now verify the equicontinuity of $G$ on $B$. Let $\sigma_i=t_n+\theta_i h_{n+1}$ and $0<\theta_2-\theta_1\leq 1$. For $y\in B$, it holds
$$
\begin{aligned}
  G(y)(\sigma_1)-G(y)(\sigma_2)  
 & =   \big(\mathrm{e}^{-\theta_1 h_{n+1}A}u_n - \mathrm{e}^{-\theta_2 h_{n+1}A}u_n\big) +    \sum_{i=1}^s\int_0^{ \theta_1 h_{n+1}} \mathrm{e}^{- (\theta_1  h_{n+1}-\sigma)A}K_y^{(i)}(\sigma) \, \mathrm{d}\sigma\\
 & \quad-  \sum_{i=1}^s\int_0^{ \theta_2 h_{n+1}} \mathrm{e}^{- (\theta_2  h_{n+1}-\sigma)A}K_y^{(i)}(\sigma) \, \mathrm{d}\sigma   \\
 & = \big(\mathrm{e}^{-\theta_1 h_{n+1}A}u_n - \mathrm{e}^{-\theta_2 h_{n+1}A}u_n\big) -\sum_{i=1}^s \int_{\theta_1 h_{n+1}}^{\theta_2 h_{n+1}} \mathrm{e}^{-(\theta_2 h_{n+1}-\sigma)A} K_y^{(i)}(\sigma)\,\mathrm{d}\sigma  \\
 &\quad +\sum_{i=1}^s\int_0^{\theta_1 h_{n+1}} \Big(\mathrm{e}^{-(\theta_2-\theta_1)h_{n+1}  A} -I\Big)\mathrm{e}^{-(\theta_1 h_{n+1} -\sigma )A} K_y^{(i)}(\sigma)\,\mathrm{d}s  =: I_1+I_2+I_3.
\end{aligned}
$$
For arbitrary $\varepsilon>0$, by the property of the strongly continuous semigroup, there exists $\delta \in (0,1)$ such that $\|I_1\|\leq \varepsilon$ for $0<\theta_2-\theta_1\leq \delta$. Let $\delta^\prime=\min\{\varepsilon,\delta\}$ and $0<\theta_2-\theta_1\leq \delta^\prime$. Using $\|\mathrm{e}^{tA}\|_{X \leftarrow X} \leq C_{\mathrm{s}}$, the term $I_2$ is bounded by
$$
\|I_2\| \leq \varepsilon s C_{\mathrm{s}}T \widetilde{K},\quad \mbox{where}~\widetilde{K}= \max_{i=1,\ldots,s}\max_{0\leq \sigma \leq h_{n+1}} \|K_y^{(i)}(\sigma)\| .
$$
Using $\|t^{\alpha}A^{\alpha}\mathrm{e}^{-tA}\|_{X\leftarrow X}\leq C_{\mathrm{s}}$ and $\|(\mathrm{e}^{-tA}-I)v\| \leq C_{\mathrm{s}} t^{\alpha}\|A^{\alpha}v\|$ for $\alpha\in(0,1)$ (see \cite[Theorem 1.4.3]{Henry81Book}), one has 
$$
\begin{aligned}
\|I_3\|&\leq C_{\mathrm{s}}\sum_{i=1}^s  \int_0^{\theta_1h_{n+1}} \big( (\theta_2-\theta_1)h_{n+1}\big)^{\frac{1}{2}} \|A^{\frac{1}{2}} \mathrm{e}^{-(\theta_1 h_{n+1}-\sigma)A}K_y^{(i)}(\sigma) \big\|\,\mathrm{d}\sigma \\
&\leq    s C_{\mathrm{s}}^2\widetilde{K}\int_0^{\theta_1h_{n+1}} \big( (\theta_2-\theta_1)h_{n+1}\big)^{\frac{1}{2}} (\theta_1 h_{n+1}-\sigma)^{-\frac{1}{2}} \,\mathrm{d}\sigma \leq 2\varepsilon^{\frac{1}{2}}  s C_{\mathrm{s}}^2 T\widetilde{K}.
\end{aligned}
$$
We conclude that the image of $B$ under $G$ is an equicontinuous family of functions. It is obvious that $G(B)$ is uniformly bounded. Therefore, we can apply the Arzela--Ascoli theorem to conclude that $G(B)$ is precompact in $B$. Finally, Schauder's fixed point theorem ensures the existence of a solution.

We proceed to establish the uniqueness of the solution. For any solution $U=G(U)$ with $U\in Y$,  using \cite[Theorem 4.3.1]{Lunardi95Book} and the fact that
$$
\begin{aligned}
\|U\|_{C([t_n,t_{n+1}];X)} &
\leq \|U-G(w_{n+1})\|_{C([t_n,t_{n+1}];X)}  + \|G(w_{n+1})\|_{C([t_n,t_{n+1}];X)} \\
& \leq 1 +  C\big( T + \|w_{n+1}\|_{C([t_n,t_{n+1}];X)}  + \|U\|_{C((-\infty,t_n];X)} + \|f(0,0,0)\|\big),
\end{aligned}
$$ 
we obtain   
$$
\begin{aligned}
\|U\|_{C^1([t_n,t_{n+1}];X)}& \leq C\big(\|A u_n\| + \|r_{n+1}\|_{C^1([t_n,t_{n+1}];X)} \big) \\
&\leq  C  \|A u_n\| + C\sum_{i=1}^s \big\|f\big( t_{ni},U_{ni},U\big(t_{ni}-\tau(t_{ni},U_{ni})\big) \big)-f(0,0,0)\big\| + C\|f(0,0,0)\| \\
&\leq C \big( \|A u_n\|  +  T + \|U\|_{C((-\infty,t_n];X)} + \|U\|_{C([t_n,t_{n+1}];X)}+\|f(0,0,0)\|\big) \\
&\leq C \big( \|A u_n\|   + \|U\|_{C((-\infty,t_n];X)} + \|w_{n+1}\|_{C([t_n,t_{n+1}];X)}+\|f(0,0,0)\| +1\big) .
\end{aligned}
$$
For any $y\in Y$, we have 
$$
\begin{aligned}
&\|G(U)-G(y)\|_{C([t_n,t_{n+1}];X)} \\
&\quad\leq  
C_{\mathrm{s}} C_f h_{n+1}\sum_{i=1}^s\Big( \|U(t_{ni})-y(t_{ni})\| +  \big\|U\big(t_{ni}-\tau(t_{ni},U(t_{ni})  )\big) -\widehat{y}\big(t_{ni}-\tau(t_{ni},y(t_{ni})  )\big) \big\|  \Big) \\
&\quad\leq C_{\mathrm{s}} C_f h_{n+1}\sum_{i=1}^s\bigg(  \|U(t_{ni})-y(t_{ni})\| +\big\|U\big(t_{ni}-\tau(t_{ni},U(t_{ni})  )\big) -U\big(t_{ni}-\tau(t_{ni},y(t_{ni})  )\big) \big\|  \\
&\quad\qquad\qquad\qquad\qquad +\big\|U\big(t_{ni}-\tau(t_{ni},y(t_{ni})  )\big) -\widehat{y}\big(t_{ni}-\tau(t_{ni},y(t_{ni})  )\big) \big\|  \bigg).
\end{aligned}
$$
Using the Lipschitz continuity of $U$, we further have
\begin{equation*}
\|G(U)-G(y)\|_{C([t_n,t_{n+1}];X)} \leq  Ch_{n+1}\|U-y\|_{C([t_n,t_{n+1}];X)}.
\end{equation*}
For sufficiently small $h_{n+1}$, the map $G$ always contracts the distance $\|U-y\|_{C([t_n,t_{n+1}];X)}$ from the solution $U$. It follows that the solution $U$ is unique, which completes the proof.
\end{proof}

Similar to the convergence analysis in Section \ref{Sec:second}, we need to consider  the numerical solution of the local problem \eqref{Eqn:secondERK-zn+1} the ERK methods of collocation type \eqref{Eqn:colERK-un+1}-\eqref{Eqn:continuous-colERK-Un+1}, which has the formula
\begin{equation}\label{Eqn:ERKcol-hatun+1}
\begin{aligned}
\begin{aligned}
\widehat{u}_{n+1} & = \mathrm{e}^{-h_{n+1}A}u_n +   h_{n+1}b_i(-h_{n+1}A) f\big(t_{ni},\widehat{U}_{ni}, u \big(t_{ni}-\tau(t_{ni},\widehat{U}_{ni}  ) \big)\big)  \\
\widehat{U}_{ni} & = \mathrm{e}^{-c_i h_{n+1}A}u_n +h_{n+1} \sum_{j=1}^s a_{ij}(-h_{n+1}A) f\big(t_{nj},\widehat{U}_{ni},u\big(t_{nj}-\tau(t_{nj},\widehat{U}_{nj}  ) \big)\big),\quad 1\leq i \leq s,
\end{aligned}
\end{aligned}
\end{equation}
and the corresponding continuous extension
\begin{equation}\label{Eqn:continuous-ERKcol-hatUn+1}
\begin{aligned}
\widehat{U}(t_n+\theta h_{n+1}) & = \mathrm{e}^{-\theta h_{n+1}A}u_n + h_{n+1}\sum_{i=1}^s b_i(\theta;-h_{n+1}A) f\big(t_{ni},\widehat{U}_{ni},u \big(t_{ni}-\tau(t_{ni},\widehat{U}_{ni}  ) \big)\big).
\end{aligned}
\end{equation}
The local error estimate is given in the next lemma.

\begin{lemma}\label{Lem:col}
Under the Assumptions \ref{Ass:sectorial}-\ref{Ass:lip}, if the function
$$
g(t)=f\big(t,u(t),u\big(t-\tau(t,u(t))\big)\big)
$$
is of class $C^{s-1,1}$ on $[t_n ,t_{n+1}]$, then the following error bounds
$$
\max_{t_n\leq t \leq t_{n+1}} \|\widehat{U}(t)-u(t)\| \leq Ch_{n+1}^{s+1},
$$
holds. The constant $C$ is independent of $h_{n+1}$. 
\end{lemma}
\begin{proof}
Let $\widetilde{e}(t)=\widehat{U}(t)-u(t)$.
Similar to the derivation in \cite[Equations (4.18)-(4.22)]{HuangOstermann25}, we can obtain that
\begin{align}
 \widehat{e}(t_n+\theta h_{n+1}) 
& =   h_{n+1} \sum_{i=1}^s b_i(\theta;- h_{n+1}A)\Big( f\big(t_{ni},\widehat{U}_{ni}, u\big(t_{ni}-\tau(t_{ni},\widehat{U}_{ni} \big)\big) - g(t_{n,i}) \Big) -\Delta_{n+1}(\theta),  \label{Eqn:haten+1}
\end{align}
where the defect $\Delta_{n+1}$ is given as
$$
\begin{aligned}
\Delta_{n+1}(\theta)= & \int_0^{\theta h_{n+1}} \mathrm{e}^{-(\theta h_{n+1}-\sigma) A} \int_0^\sigma \frac{(\sigma-\xi)^{s-1}}{(s-1)!} g^{(s)}(t_n+\xi) \,\mathrm{d} \xi 
\,\mathrm{d} \sigma  \notag \\
& - h_{n+1} \sum_{i=1}^s b_i(\theta;- h_{n+1} A) \int_0^{c_i h_{n+1}} \frac{(c_i h_{n+1}-\sigma)^{s-1}}{(s-1)!} g^{(s)}(t_n+\sigma) \,\mathrm{d} \sigma. 
\end{aligned}
$$
Therefore, we have
$$
\max_{0\leq \theta \leq 1}\|\Delta_{n+1}(\theta)\| \leq Ch_{n+1}^{s+1}.
$$

Using \eqref{Eqn:boundedness-bitheta} and the Lipschitz continuity of $f$, $u$ and $\tau$, we obtain from \eqref{Eqn:haten+1},
$$
\begin{aligned}
&\max_{0\leq \theta \leq 1}\|\widetilde{e}(t_n+\theta h_{n+1})\| \\
 & \qquad\leq Ch_{n+1}\sum_{i=1}^s\Big( \|\widetilde{e}(t_{ni})\|+ \big\|  u\big(t_{ni}-\tau(t_{ni},\widehat{U}_{ni} \big)\big) - u\big(t_{ni}-\tau(t_{ni},u(t_{ni}) \big) \big\| \Big) +  Ch_{n+1}^{s+1} \\
&\qquad \leq Ch_{n+1} \max_{0\leq \theta \leq 1}\|\widetilde{e}(t_n+\theta h_{n+1})\| +  Ch_{n+1}^{s+1},
\end{aligned}
$$
Thus, for $h_{n+1}$ suﬀiciently small, we obtain
$$
\max_{0\leq \theta \leq 1}\|\widetilde{e}(t_n+\theta h_{n+1})\|  \leq Ch_{n+1}^{s+1},
$$
which completes the proof.
\end{proof}

The convergence result of the ERK methods of collocation type \eqref{Eqn:colERK-un+1}-\eqref{Eqn:continuous-colERK-Un+1} is stated below. Its proof is similar to that of Theorem \ref{Thm:second-order}. For the sake of brevity, we omit the details here.

\begin{theorem}\label{Thm:col}
Under the Assumptions \ref{Ass:sectorial}-\ref{Ass:lip}, let $g$ be of class $C^{s-1,1}$ on the intervals $[t_j,t_{j+1}]$, $j=0,\ldots,N-1$. Consider for the numerical solution of the initial value problem \eqref{Eqn:problem} an ERK method of collocation type \eqref{Eqn:colERK-un+1}-\eqref{Eqn:continuous-colERK-Un+1}. Then for suﬀiciently small $h$, the error bound
$$
\|u_n-u(t_n)\|\leq Ch^s 
$$
holds uniformly on $0\leq t \leq T$. The constant $C$ depends on $T$, but is independent of the step size sequence.
\end{theorem}

Provided that the underlying quadrature rule is of order $s+1$, i.e.,
$$
\sum_{i=1}^s b_i(0)c_i^s = \frac{1}{s+1},
$$
the ERK methods of collocation type \eqref{Eqn:colERK-un+1}-\eqref{Eqn:continuous-colERK-Un+1} can achieve order $s+1$. This superconvergence is stated as below. Its proof is quite similar to that of  \cite[Theorem 5.1]{HuangOstermann25}.

\begin{theorem}\label{Thm:col}
Under the Assumptions \ref{Ass:sectorial}-\ref{Ass:lip}, let $g$ be of class $C^{s,1}$ on the intervals $[t_j,t_{j+1}]$, $j=0,\ldots,N-1$. Consider for the numerical solution of the initial value problem \eqref{Eqn:problem} the ERK methods of collocation type \eqref{Eqn:colERK-un+1}-\eqref{Eqn:continuous-colERK-Un+1} whose underlying quadrature rule is of order $s + 1$. Let the step size sequence $\{h_j\}_{j=1}^N$ satisfy the condition $h_j\leq \varrho h_{j+1}$ with $\varrho >1$ for all $j$.
Then for suﬀiciently small $h$, the error bound satisfies
$$
\|u_n-u(t_n)\|\leq CC_{\mathrm{S}}h^{s+1},
$$
uniformly on $0\leq t \leq T$. In general, the size of $C_{\mathrm{S}}$ depends on the chosen step size sequence as follows
$$
1\leq C_{\mathrm{S}} \leq \varrho \ln \frac{T}{\min_{1\leq j \leq N} h_j}+2.
$$
However, when the step sizes are constant or when the operator $A$ and the space $X$ satisfy certain conditions (see \cite[Remark 1]{HuangOstermann25}), $C_{\mathrm{S}}$ is independent of the step size sequence. On the other hand, the constant $C$ depends on $T$, but not on the step size sequence.
\end{theorem}

\section{Numerical experiments and implementation}\label{Sec:experiments}
In this section, we first comment on the implementations of ERK methods. Then some numerical experiments are presented to illustrate the convergence results obtained in the previous sections.

\subsection{Implementation issues}\label{Sec:issues}
As mentioned before, although the underlying method \cite[Equation (5.3)]{HochbruckOstermann05:1069} is explicit, the second order ERK method \eqref{Eqn:secondERK-un+1}-\eqref{Eqn:continuous-secondERK-Un+1} is implicit when overlapping occurs. It is common to determine the continuous extension of the solution by iteration using a predictor-corrector method (cf.~\cite{HartungKrisztin06Book}). Recall the notation introduced in \eqref{Eqn:b1b2}
$$
b_1(\theta;-h_{n+1}A)=\theta\varphi_1(- \theta h_{n+1} A)-\tfrac{1}{c_2}\theta^2\varphi_2(-\theta h_{n+1}A) , \quad
b_2(\theta;-h_{n+1}A)=\tfrac{1}{c_2} \theta^2\varphi_2(-\theta h_{n+1}A).
$$
The following pseudo-code performs one step in the predictor-corrector mode (with $m$ corrections). In practice, the number $m$ of corrections need not be fixed a prior; one can stop when the difference between two successive $\widehat{u}_{n2}$ falls below a prescribed tolerance.

\begin{algorithm}[H]
\caption{Predictor-Corrector$^m$ Mode for \eqref{Eqn:secondERK-un+1}-\eqref{Eqn:continuous-secondERK-Un+1}}
\begin{algorithmic}[0]
\State \textbf{Step 1:} Predictor
\State $G_{n1}=f (t_n,u_n,U  (t_n-\tau(t_n,u_n  )  ) )$
\State $U_{n2}=\mathrm{e}^{-c_2 h_{n+1}A}u_n + c_2  h_{n+1}\varphi_1(-  c_2 h_{n+1} A) G_{n1}$
\If{$t_{n2}-\tau(t_{n2},U_{n2}) \leq t_n$}
    \State $G_{n2}=f(t_{n2},U_{n2},U(t_{n2}-\tau(t_{n2},U_{n2})))$
    \Else
    \State $G_{n2}=f(t_{n2},U_{n2},u_n)$    
\EndIf
\State \textbf{Step 2:} Correction by iteration is needed if $t_{n2}-\tau(t_{n2},U_{n2}) > t_n$
\State{$\theta_2=\frac{t_{n2}-\tau(t_{n2},U_{n2}) - t_n}{h_{n+1}}$}
\For{$r=1,\ldots,m$}
\State{$\widehat{u}_{n2}= \mathrm{e}^{-\theta_2 h_{n+1}A}u_n  + h_{n+1}b_1(\theta_2;-h_{n+1}A)  G_{n1} + h_{n+1}  b_2(\theta_2;-h_{n+1}A) G_{n2}$}
\State{$G_{n2}=f(t_{n2},U_{n2},\widehat{u}_{n2})$}
\EndFor
\State \textbf{Step 3:} Computation of the continuous extension to $[t_n,t_{n+1}]$
\State{$U(t_n+\theta h_{n+1})  = \mathrm{e}^{-\theta h_{n+1}A}u_n + b_1(\theta;-h_{n+1}A) G_{n1}+ h_{n+1}b_2(\theta;-h_{n+1}A)G_{n2}$}
\end{algorithmic}
\end{algorithm}

Since the ERK methods of collocation type (with $s\geq 2$) are implicit for standard semilinear parabolic problems, their implementation is typically more involved. We additionally need to conduct a fixed point iteration to evaluate the value of $U(t_{ni}-\tau(t_{ni},U(t_{ni}))$ if $t_{ni}-\tau(t_{ni},U(t_{ni})) >t_n$. The following pseudo-code performs one step in the predictor-(evaluation-corrector)$^m$ mode.

\begin{algorithm}[H]
\caption{Predictor-(Evaluation-Corrector)$^m$ Mode for \eqref{Eqn:colERK-un+1}-\eqref{Eqn:continuous-colERK-Un+1}}
\begin{algorithmic}[0]
\State \textbf{Step 1:} Predictor
\State Set $U_{ni}^{(0)}=u_n$ for $i  =  1,\ldots,s$
\State \textbf{Step 2:}  Evaluation-Corrector
\For{$r=1,\ldots,m$}
    \State{$\bullet$~Evaluation:}
    \State Set $Y=\emptyset$
    \For{$i=1,\ldots,s$}
    \State{$s_i=t_{ni}-\tau(t_{ni},U_{ni}^{(r-1)})$}
    \If{$s_i \leq t_n$}
    \State{$X_i=U(s_i)$} 
    \Else
    \State{$\theta_i=\tfrac{s_i-t_n}{h_{n+1}}$}
    \State{$Y=Y\cup \{i\}$}
    \EndIf
    \EndFor
    \If{$Y\ne \emptyset$}
    \State{Solve $X_i=\mathrm{e}^{-\theta_ih_{n+1}A}u_n+\sum_{j=1}^sb_j(\theta_i;-h_{n+1}A)f(t_{nj},U_{nj}^{(r-1)},X_j),\quad \mbox{for}~i\in   Y$}
    \EndIf    
    \State{$\bullet$~Correction: $U_{ni}^{(r)}=\mathrm{e}^{-c_ih_{n+1}A}u_n+\sum_{i=1}^sb_j(c_i;-h_{n+1}A)f(t_{ni},U_{ni}^{(r-1)},X_i),\quad \mbox{for}~i=1,\ldots,s$}
\EndFor
\State \textbf{Step 3:} Computation of the continuous extension to $[t_n,t_{n+1}]$
\State{$U(t_n+\theta h_{n+1})  = \mathrm{e}^{-\theta h_{n+1}A}u_n + \sum_{i=1}^s b_i(\theta;-h_{n+1}A) f(t_{ni},U_{ni}^{(m)},X_i)$}
\end{algorithmic}
\end{algorithm}

The convergence result for high order methods in Theorem \ref{Thm:col} requires that $g(t)=f\big(t,u(t),u\big(t-\tau(t,u(t))\big)\big)$ is sufficiently smooth on each interval $[t_j,t_{j+1}]$. However, this composition generally exhibits low regularity at certain points, due to the fact that the solution $u(t)$ does not connect smoothly to the initial function $\phi(t)$ (see Remark \ref{Rem:second}). In practice, the low regularity points ought to be included in the mesh to avoid the loss of accuracy. Consider the following spatial discretization system of problem \eqref{Eqn:problem}, arising for instance from finite difference or finite element methods:
\begin{equation*}\label{Eqn:spatial-problem}
\left\{\begin{aligned}
&\mathbf{U}^\prime(t) + \mathbf{A} \mathbf{U}(t) = \mathbf{f}\big(t,\mathbf{U}(t),\mathbf{U}   \big(t-\widetilde{\tau}(t,\mathbf{U}(t))\big) \big), && 0\leq t \leq T, \\
& \mathbf{U}(t) =\bm{\Phi}(t) &&t\leq 0, \end{aligned}\right.
\end{equation*}
where $\mathbf{U}(t) \in \mathbb{R}^m$ is the approximation of the solution $u(t) \in X$. The nonlinearity $\mathbf{f} : [0,T] \times \mathbb{R}^m \times \mathbb{R}^m \to \mathbb{R}^m$ and the delay $\widetilde{\tau}: [0,T] \times \mathbb{R}^m \to \mathbb{R}_{\geq 0}$ are obtained via spatial discretization of $f$ and $\tau$, respectively. The matrix $\mathbf{A} \in \mathbb{R}^{m \times m}$ is the discretization of a differential operator. This leads to a stiff system of state-dependent delay differential equations.
If $u'(0^-) \ne \phi'(0^+)$, then a consistent semi-discrete solution of \eqref{Eqn:problem} reproduces this lack of smoothness, that is, $\mathbf{U}'(0^-) \ne \bm{\Phi}'(0^+)$, where $\bm{\Phi}$ is the spatial discretization of the initial data $\Phi$. As is well known \cite{BellenMaset09:1,FeldsteinNeves84:844}, this derivative jump at $t = 0$ is propagated and ``smoothed'' by the lag term $t - \widetilde{\tau}(t, \mathbf{U}(t))$. There are only finitely many critical points. We label these points as an increasing sequence $0=\xi_0<\xi_1 < \xi_2 <\cdots < \xi_{\ell} \leq T$. Each discontinuity point $\xi_j$ ($j\neq 0$) is a descendent of some previous point $\xi_i$, satisfying the relation
$$
\xi_j-\widetilde{\tau}(\xi_j,\mathbf{U}(\xi_j)) = \xi_i,\quad 0\leq i < j \leq \ell.
$$
The locations of these points cannot be computed a prior since their unknown locations $\xi_j$ depend implicitly on the also unknown solution $\mathbf{U}$. Extensive work has been devoted to tracking discontinuities in state-dependent DDEs; see \cite{BellenMaset09:1,BellenZennaro13Book,XuHuang19:314} and the references therein. We consider the switching function method developed in \cite{{KarouiVaillancourt94:37}}. Suppose that the steps $\mathbf{Y}_1,\ldots,\mathbf{Y}_n$ were aleardy obtained by an ERK method of order $p$, and the approximate discontinuity points found so far are $0=\widetilde{\xi}_0<\widetilde{\xi}_1 < \widetilde{\xi}_2 <\cdots < \widetilde{\xi}_{\vartheta} < T$.
\begin{description}
\item[Step 1:] Compute the next approximate value $\mathbf{Y}_{n+1}$ ($\approx \mathbf{U}(t_{n+1})$) using the ERK method with a given step size $h_{n+1}$.
\item[Step 2:] For $i=1,\ldots,\vartheta$ find some $i$ such that
$$
\big( t_n-\widetilde{\tau}(t_n,\mathbf{Y}_n) -\widetilde{\xi}_i \big) \big(  t_{n+1}-\widetilde{\tau}(t_{n+1},\mathbf{Y}_{n+1}) -\widetilde{\xi}_i\big) <0.
$$
If such $i$ does not exist, then the current step size $h_{n+1}$ and solution $\mathbf{Y}_{n+1}$ are accepted and the algorithm proceeds to the next integration step. Otherwise, we proceed with Step 3.
\item[Step 3:] Construct an interpolation polynomial $Q(t)$ of degree $p-1$ satisfying 
$$
Q(t_k)=t_k-\widetilde{\tau}(t_k,\mathbf{Y}_k) -\widetilde{\xi}_i,\quad k=n-p+1,\ldots,n. 
$$
Use the bisection method to find the root $\widetilde{\xi}$ of $Q(t)$ in the interval $(t_n,t_n+h_{n+1})$, and set $\widetilde{\xi}_{\vartheta +1 } = \widetilde{\xi}$.
\item[Step 4:] Set $t_{n+1}=\widetilde{\xi}$ as the next mesh point and compute the corresponding solution $\mathbf{Y}_{n+1}$.
\end{description}

\subsection{Convergence tests}
We test the convergence rates of various ERK methods developed in previous sections. 
The first order method refers to the exponential Euler method \eqref{Eqn:pre-euler-un+1}-\eqref{Eqn:continuous-euler-Un+1}, while the second order method corresponds to the method given in \eqref{Eqn:secondERK-un+1}-\eqref{Eqn:continuous-secondERK-Un+1} with $c_2=1$.
The third and fourth order methods are of collocation type \eqref{Eqn:colERK-un+1}-\eqref{Eqn:continuous-colERK-Un+1}, with the third order method using collocations points $c_1=1/3,c_2=2/3,c_3=1$, and the fourth order method employing Gauss--Lobatto collocation points $c_1=0,c_2=1/2,c_3=1$. 

\begin{example}\rm
We begin by investigating the following one-dimensional parabolic problem with
known exact solution
\begin{equation}\label{Eqn:1}
\partial_t u - \partial_{xx} u = \frac{1}{1+u^2 + \big( u (t-  \tau(t,u)   ) \big)^2} + \Psi(x,t)\quad \mbox{with}~\tau(t,u)=(1-t)\|u\|_{L^2}^2
\end{equation}
for $u=u(t,x)$, where $t\in[0,1]$ and $x\in [0,1]$, subject to the homogeneous Dirichlet boundary conditions. The source function $\Psi$ is determined by the exact solution of the problem 
$$
u(t,x)=\mathrm{e}^{t}x(1-x),\quad t \in [-\tfrac{1}{30},1].
$$ 
\end{example}

\begin{figure}[H]
    \centering
    \includegraphics[width=0.5\textwidth]{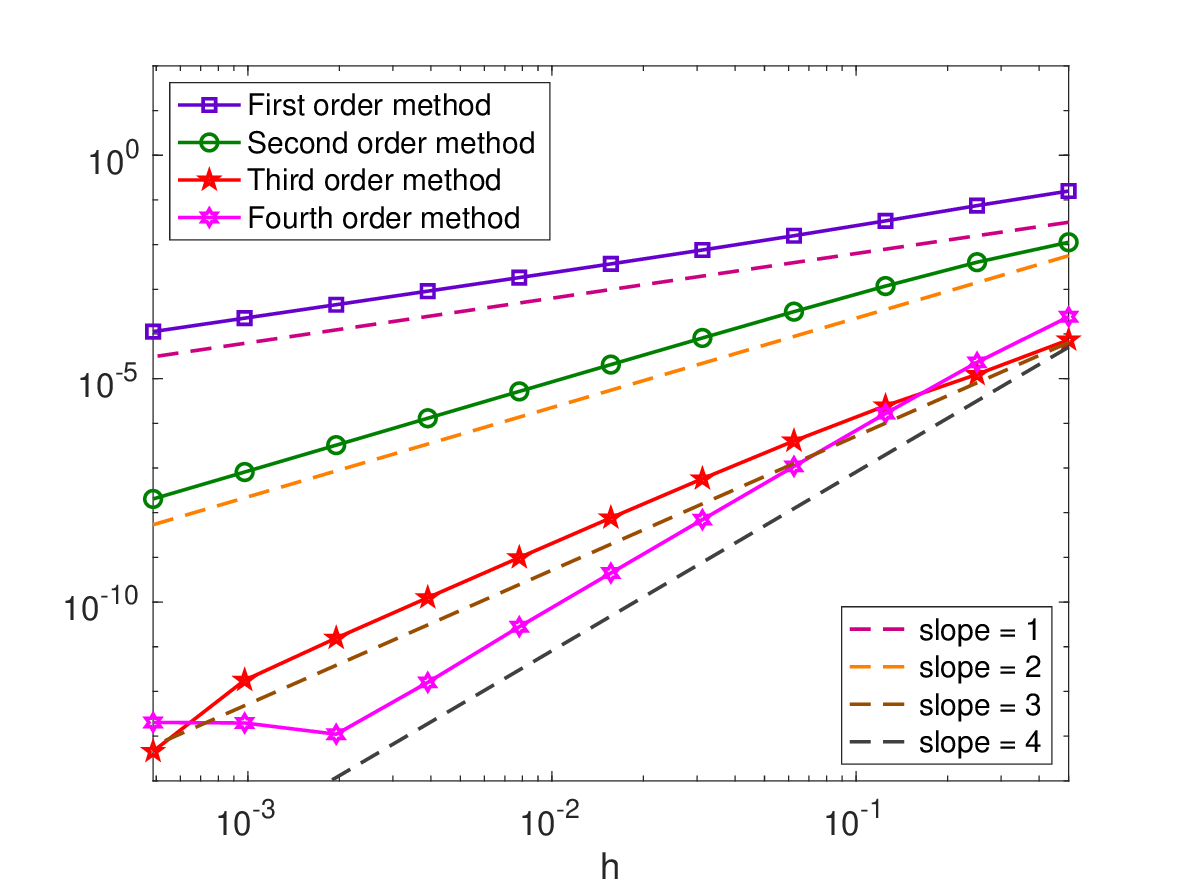}
    \caption{The convergence rates of ERK methods for \eqref{Eqn:1}. The errors are measured at $T=1$ in the $L^2(\Omega)$ norm.}
    \label{Fig:1} 
\end{figure}

We apply a standard finite difference method with $n = 200$ grid points to discretize  the problem in space. The resulting products of matrix functions with vectors are computed by the fast Fourier transform. The convergence rates of the ERK methods are presented in Figure~\ref{Fig:1}. The observed orders evidently are in line with our theoretical analysis.

\begin{example}\rm
In this example, we consider the following problem 
\begin{equation}\label{Eqn:2}
\partial_t u - \partial_{xx} u = \frac{1}{1+u^2 + \big( u (t-  \tau(t,u)   ) \big)^2} \quad \mbox{with}~\tau(t,u)=t- \frac{0.9 t}{1+\|u\|_{L^2}^2}
\end{equation}
for $u=u(t,x)$, where $t\in [0,1]$ and $x\in [0,1]$, subject to the homogeneous Dirichlet boundary conditions. Note that the delayed argument $t-\tau(t,u)$ is non-negative and the delay vanishes at $t=0$. The initial condition is given by $\phi(x)=x(1-x)$.
\end{example}

\begin{figure}[H]
    \centering
    \includegraphics[width=0.5\textwidth]{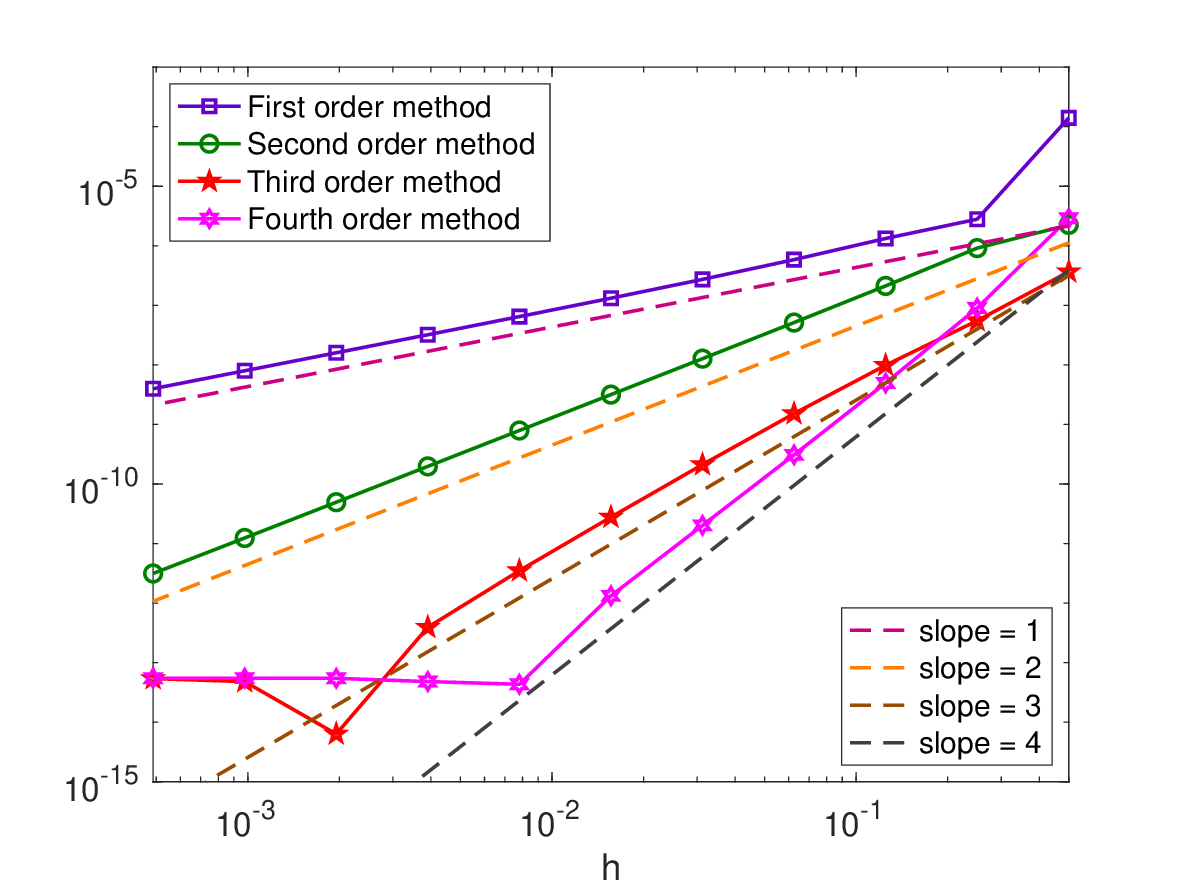}
    \caption{The convergence rates of ERK methods for \eqref{Eqn:2}. The errors are measured at $T=1$ in the $L^2(\Omega)$ norm.}
    \label{Fig:2} 
\end{figure}

We apply a standard finite difference method with $n = 200$ grid points to discretize  the problem in space. In this example the exact solution is unknown.  The reference solution is computed by the ERK method of Gauss collocation type using the constant step size $h=2^{-18}$. The errors of the ERK methods in this example are presented in Figure~\ref{Fig:2}. The numerical results clearly exhibit the expected convergence rates.

\begin{example}\rm
In the last example, we consider the following problem 
\begin{equation}\label{Eqn:3}
\partial_t u - \partial_{xx} u = \frac{1}{1+u^2 + \big( u (t-  \tau(t,u)   ) \big)^2} \quad \mbox{with}~\tau(t,u)= \frac{2}{3+\|u\|_{L^2}^2}
\end{equation}
for $u=u(t,x)$, where $t\in [0,1]$ and $x\in [0,1]$, subject to the homogeneous Dirichlet boundary conditions. The initial condition is given by $\phi(t,x)=\mathrm{e}^t x(1-x)$ with $t\in[-1,0]$. Note that the delay does not vanish at $t=0$. 
\end{example}

We apply a standard finite difference method with $n = 200$ grid points to discretize  the problem in space, leading to a stiff system of state-dependent DDEs.  Since the delay does not vanish at $t=0$, potential derivative discontinuities must be tracked and incorporated into the time mesh. The reference solution is computed using an ERK method based on Gauss collocation with a default time step size $h = 2^{-19}$. To capture potential discontinuities, we employ the switch function method to adaptively adjust the time step. As a result, a discontinuity is detected at $t = 0.664973949550472$.

\begin{figure}[H]
    \centering
    \includegraphics[width=0.49\textwidth]{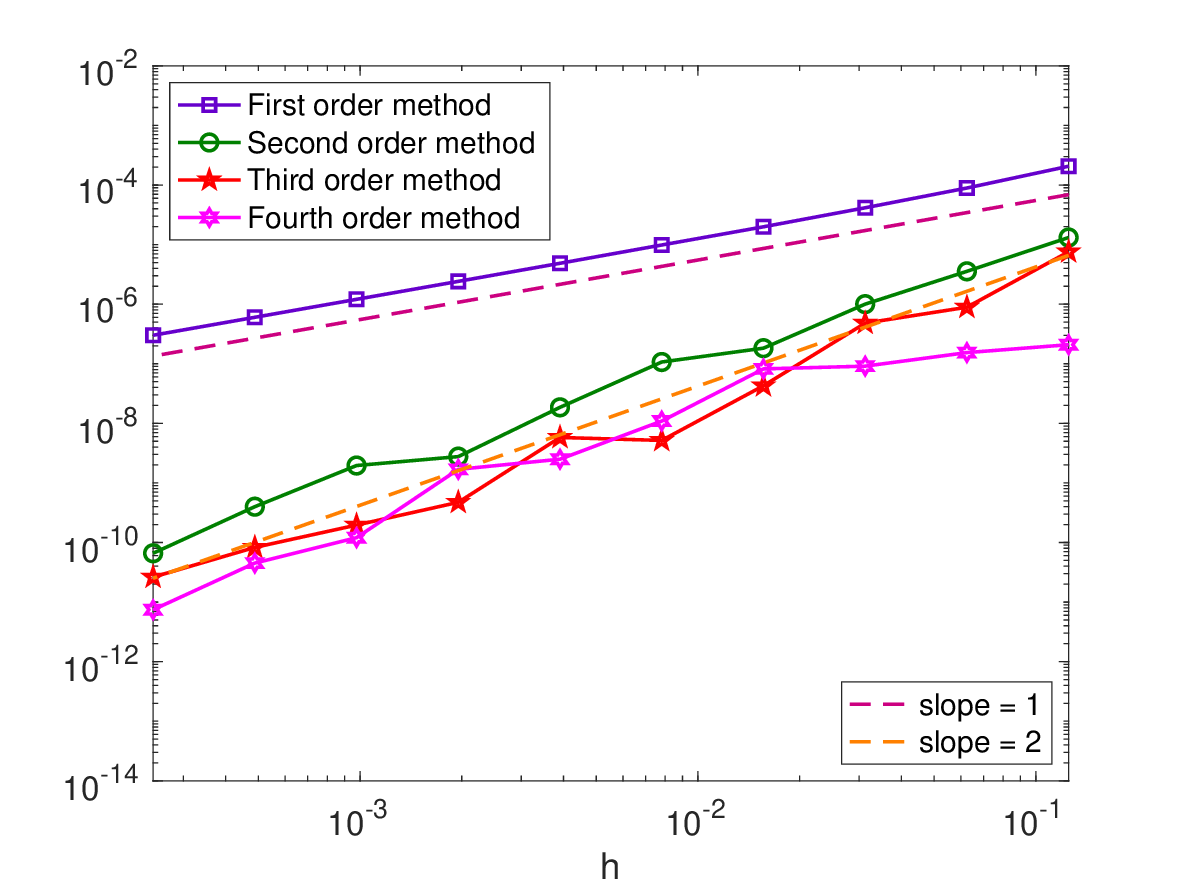}
    \includegraphics[width=0.49\textwidth]{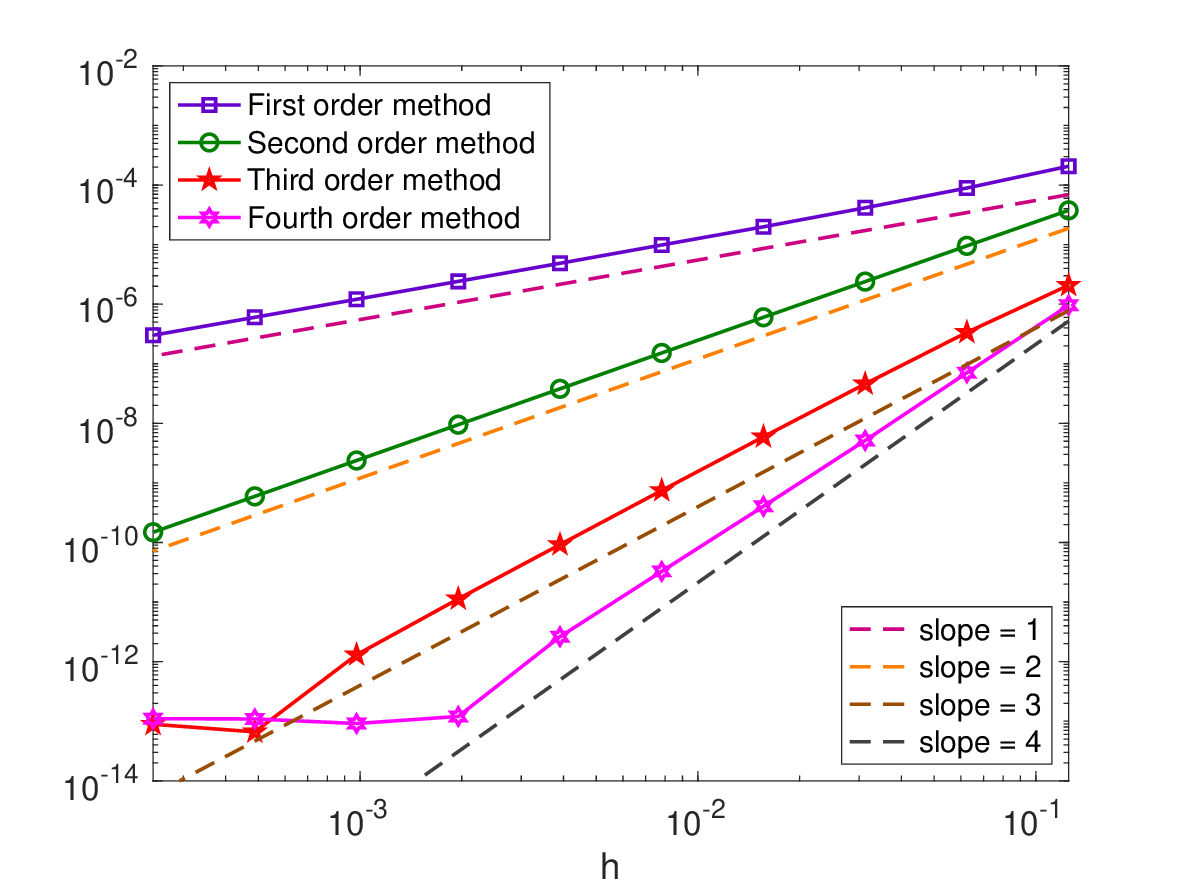}
    \caption{The convergence rates of ERK methods for \eqref{Eqn:3}. The errors are measured at $T=1$ in the $L^2(\Omega)$ norm. Left: fixed step size. Right: the step size is adjusted via the switch function method to capture the discontinuity, which is only applied for methods of order at least two.
}
    \label{Fig:3} 
\end{figure}

We investigate the convergence behavior of ERK methods under two scenarios:
(i) using a constant step size without capturing the discontinuity, and
(ii) using the same step size by default, but locally adjusting it based on the switch function method when a discontinuity is detected. The convergence rates evaluated at $T=1$ in the $L^2(\Omega)$ norm are presented in Figure~\ref{Fig:3}. It is observed that the convergence rate is significantly reduced when the discontinuity is not captured by the time mesh, with at most second order convergence being observed. In contrast, incorporating the discontinuity into the mesh enables the ERK methods to achieve the expected order of convergence.

\section*{Acknowlegments}
Qiumei Huang is supported by the National Natural Science Foundation of China (No.~12371385). 
Gangfan Zhong is supported by the China Scholarship Council (CSC) joint Ph.D. student
scholarship (Grant 202406540082).

\bibliographystyle{abbrv}

\bibliography{/Volumes/Bibliography/gfzhong_bib/Mathematics.bib}


\end{document}